\documentclass[11pt]{article}   % For Latex2e
\usepackage{amssymb,amscd,latexsym} 
\usepackage{graphicx,color}  % For Latex2e
\usepackage[pagebackref]{hyperref}
\usepackage{amsmath}
\usepackage{amsthm}
\usepackage{mathdots}
\usepackage[all]{xy}
\usepackage{url}
\newcommand{\llar}{-\kern-5pt-\kern-5pt\longrightarrow}
\newcommand{\surjects}{\twoheadrightarrow}

\newtheorem{Theorem}{Theorem}[section]
\newtheorem{Lemma}[Theorem]{Lemma}

\newtheorem{Proposition}[Theorem]{Proposition}
\newtheorem{Remark}[Theorem]{Remark}
\newtheorem{Example}[Theorem]{Example}

\newtheorem{Definition}[Theorem]{Definition}

\def\sqr#1#2{{\vcenter{\hrule height.#2pt
        \hbox{\vrule width.#2pt height#1pt \kern#1pt
            \vrule width.#2pt}
        \hrule height.#2pt}}}
\def\phi{\varphi}
\def\demo{\noindent{\bf Proof. }}
\def\square{\mathchoice\sqr64\sqr64\sqr{4}3\sqr{3}3}
\def\qed{\hspace*{\fill} $\square$}

\def\xx{{\bf x}}
\def\yy{{\bf y}}

\def\XX{{\bf X}}

\def\vv{{\bf v}}

\def\codim{{\rm codim}\,}

\def\rk{{\rm rank}\,}

\def\restr{{\kern-1pt\restriction\kern-1pt}}

\def\pp{{\mathbb P}}

\begin{document}
%\begin{titlepage}
\begin{center}
{\Large{\bf\sc Symmetry preserving degenerations of the generic symmetric matrix}}
\footnotetext{AMS Mathematics
	Subject Classification (2010   Revision). Primary 13C40, 13D02, 14E05, 14M12; Secondary  13H10,14M05.}

\vspace{0.3in}

{\large\sc Rainelly Cunha}\footnote{ Under a CAPES-PNPD post-doctoral fellowship (1723357/2017) from the Federal University of Sergipe.} \quad\quad
{\large\sc Zaqueu Ramos}%\footnote{Partially supported by a CNPq post-doctoral fellowship (151229/2014-7)}
 \quad\quad
 {\large\sc Aron  Simis}\footnote{Partially supported by a CNPq grant (302298/2014-2). Parts of this work were done while this author held a FAPESP Senior Visiting Research Scholarship at ICMC-USP (S\~ao Carlos, SP), for which he is grateful.}

\end{center}
%\end{titlepage}

%\tableofcontents
%\newpage

\begin{abstract}

One considers certain degenerations of the generic symmetric matrix over a field $k$ of characteristic zero and the main structures related to the determinant $f$ of the matrix, such as the ideal generated by its partial derivatives, the polar map defined by these derivatives and its image $V(f)$, the Hessian matrix, the ideal and the map given by the cofactors, and the dual variety of $V(f)$. 
%Cases where the degenerated determinant has non-vanishing Hessian determinant show that the former is a factor of the latter with the (Segre) expected multiplicity, a result treated by Landsberg-Manivel-Ressayre by geometric means.
%Another byproduct is an affirmative answer to a question of F. Russo concerning the codimension in the polar image of the dual variety to a hypersurface.
A complete description of these structures is obtained.
\end{abstract}

\section*{Introduction}

Let $m\geq 3$ denote an integer and $k$ a field. We consider the $m\times m$ generic symmetric matrix:
\begin{equation}\label{symgeneric}
\mathcal{S}:=\left(
\begin{matrix}
x_{1,1}&x_{1,2}&\ldots & x_{1,m-1} & x_{1,m}\\
x_{1,2}&x_{2,2}&\ldots & x_{2,m-1} & x_{2,m}\\
\vdots &\vdots &\vdots &\vdots & \vdots\\
x_{1,m-1} & x_{2,m-1}&\ldots & x_{m-1,m-1} & x_{m-1,m}\\
x_{1,m} & x_{2,m}&\ldots & x_{m-1,m} & x_{m,m}
\end{matrix}
\right),
\end{equation}
where the entries are indeterminates over $k$.
Let $S:=k[x_{i,j}\,|\, 1\leq i\leq j \leq m]$ denote the polynomial ring over $k$ generated by the entries of $\mathcal{S}$.
By a ``degeneration'' of  $\mathcal{S}$ we mean the matrix $\mathcal{DS}$ obtained by applying to the entries of $\mathcal{S}$ a given $k$-endomorphism of $S$.
Thus, e.g., the matrix
$$
\left(
\begin{matrix}
a_1 & 0\\
a_2 & a_3
\end{matrix}
\right),
$$
with $a_2\neq 0$, is not a degeneration of the $2\times 2$ generic symmetric matrix no matter what values the entries have -- although it is a perfectly acceptable degeneration of the $2\times 2$ generic matrix.
Indeed, a degeneration $\mathcal{DS}$ in the present sense will preserve the symmetric nature of the original  $\mathcal{S}$, therefore establishing some a priori constraints.
In addition, the $k$-endomorphisms considered here will be {\em coordinate-like}, i.e., induced by a $k$-linear map of the $k$-vector space $S_1$ that maps a variable to another variable or to $0$.
Throughout we will often make the abuse of thinking of a degeneration as being either the action by the $k$-endomorphism or the resulting matrix itself.

Degenerations as above of the generic matrix have been considered by Merle and Giusti \cite{Merle} and by Eisenbud in \cite{Eisenbud2}.
The sort of degeneration we consider here has also been dealt with in \cite{LiZaSi} for the case of the generic square matrix. 
As it will become clear along the development, there are some genuine differences that require appropriate changes in the symmetric environment, not to mention the (at least psychological) discomfort of being no longer free to identify entries with variables in a bijective way.
Besides, the numerical invariants in the two situations will often diverge as one naturally expects. 
For this reason we will refrain from any blind reference to the arguments in \cite{LiZaSi} and instead reinstate proofs {\em ab initio} whenever required.

A major question is when the Hessian of the degeneration $f:=\det\mathcal{DS}$ vanishes. The general question of the vanishing of a hyperurface has a venerable history since the days of Hesse (\cite{Hesse1}, \cite{Hesse2}) and Gordan--Noether (\cite{GN}), subsequently studied by several mathematicians of the Italian school (\cite{Per}, \cite{Fran1}, \cite{Fran2}, \cite{Perm1}, \cite{Perm2}) and more recently by C. Ciliberto, R. Gondim, F. Russo, G. Staglian\`o (\cite{Cili}, \cite{CRS}, \cite{GonRuss}, \cite{GRS}).
In this paper the focus is on the class of determinantal hypersufaces arising from degenerations of the generic symmetric matrix.
The spirit stays closer to \cite{CRS}, with a mixed study of vanishing and non-vanishing situations, where the non-vanishing case favors a homaloidal phenomenon, while the vanishing one leads to a deep discussion of geometric invariants, including the structure and dimension of the dual variety.

Perhaps a novel point here as compared to the typical Hessian literature is a detailed consideration of the {\em gradient ideal} $J$ of $f$ (i.e., the ideal generated by the partial derivatives of $f$) and its close relationship to the ideal generated by the cofactors of the entries of $\mathcal{DS}$, as opposed to merely looking at the polar map.

We will consider two basic models of a matrix degeneration, one for the non-vanishing Hessian (in fact, homaloidal) case and one for the vanishing case.
Although the two cases are as crudely apart as they could be, the difference between the structure of the corresponding dual varieties is subtler as will be seen. 

For the first one we take the simplest coordinate-like degeneration of (\ref{symgeneric}), obtained via an endomorphism that fixes all variables appearing as entries, except one, and maps the latter to one of the other variables. 
Although  such an endomorphism has as obvious kernel the principal ideal generated by a difference of two variables, the resulting degenerations may not all share the same algebraic and homological properties of the kind we wish to consider. The relative position of the two variables may possibly inflict quite a bit of diversity on some of those properties  (see Example~\ref{bad_cloning}).
This is in striking difference with the fully generic square matrix, where the relative position of the two slots is irrelevant as long as they are not on the same row or column.
For example, one could pick an entry $x_{i,j}$ off the main diagonal of $\mathcal S$, such that $i+j$ is even,  and map it onto the entry on the $(\frac{i+j}{2}, \frac{i+j}{2})$ slot on the main diagonal, while fixing the remaining variables. Repeating this procedure along each anti-diagonal will eventually land us on a Hankel matrix.
This procedure entangles a totally different situation, which we hope to consider  in a future work.

The second model of degeneration will fix a subset of the variables while mapping the remaining ones to zero. This degeneration by coordinate hyperplanes has been considered by Merle and Giusti. However, by and large,  some of the algebraic properties may behave quite unexpectedly depending on the configuration affected by the map - a phenomenon already found in the first case above.
Thus, there is a certain strategy in the choice of the slots. Of course, it all depends on what sort of algebraic or geometric invariants one wants to analyze in the degeneration. 
%For lack of better terminology, we refer informally to such a degeneration as a {\em strategic sparsing}. 

In both situations of the Hessian status the overall goal is to understand the nature of algebraic or geometric gadgets commonly attached to the matrix  degeneration $\mathcal{DS}$.
Along with the structure of the Hessian $h(f)$ of $f=\det\mathcal{DS}$, as mentioned above we study the corresponding gradient ideal $J$, the ideal of cofactors (submaximal minors) $P$, the polar map defined by the partial derivatives of $f$ and its image $V(f)$, and the dual variety of $V(f)$.
Pretty much in general, one draws quite a bit upon the properties of the inclusion $J\subset P$ (see Proposition~\ref{GolMar}). The notation $P$ for the ideal of the cofactors expresses the wishful expectation that it ought to be prime in many cases. A characteristic behavior of $P$ is that it has maximal analytic spread, a property stated in \cite[Lemma 3.3]{LiZaSi}, as borrowed from Bruns--Vetter book (\cite[Theorem 10.16 (b)]{B-V}).

We now briefly describe the results obtained.

Section~\ref{basic} is a short piece establishing the algebraic terminology and notation.

Section~\ref{CLONING} deals with the simplest sort of coordinate-like degeneration (cloning). It is divided into two subsections. The main goal of the first subsection is to establish that the polar map associated to the determinant $f$ is birational (i.e., Cremona).
As a vehicle towards this result we prove two fundamental facts: (1)  the Hessian $h(f)$ of $f$ does not vanish; (2)  the linear syzygies of the gradient ideal $J$ of $f$ have maximal rank. 
The first of these facts is tantamount to the property that the ideal $J$ has maximal analytic spread or, equivalently, that the partial derivatives of $f$ are algebraically independent over $k$. Unfortunately, it is not clear how to establish any of these at the outset. Instead we resort to a method of specializing the Hessian matrix to a block-diagonal matrix where the blocks have well-known non-vanishing determinants.
As for the rank of the linear syzygies of $J$, we know that it is at most ${{m+1}\choose 2}-2$ (the rank of the entire syzygy matrix). Here, in contrast to the generic symmetric matrix, there will many minimal syzygies of degree $2$ as well. Fortunately, the Cauchy identity involving the matrix of cofactors affords as many linear relations as required upon identifying cofactors with partial derivatives. The problem remains to prove that the linear syzygies obtained in this manner are indeed independent, which requires a slight tour de force.

The second subsection deals with the ideal  $P\subset R$ of submaximal minors.
We first prove that this ideal is prime by proving, more strongly, that the residue ring $R/P$ is normal.
In the fully generic case considered in \cite{LiZaSi} the primeness of the ideal of submaximal minors of the cloned matrix used a result of Eisenbud drawn upon the $2$-generic property (\cite{Eisenbud2}).  Since the generic symmetric matrix is not $2$-generic, we were forced to devise an alternative to prove normality.
Next, we show that $P$ is the minimal component in a primary decomposition of the gradient ideal $J$; more precisely, $J$ is a double structure on the variety defined by $P$, with a unique embedded component of codimension $2(m-1)$ supported on a linear space.
Algebraically, this is quite a common situation where one can ask whether $J$ is actually a reduction of the prime ideal $P$ (see, e.g., \cite[Conjecture 3.16 and Corollary 3.17 (iii)]{MAron}).
The answer is negative and in order to prove this we first show that the image of the (birational) map defined by the cofactors is a hypersurface of degree $m-1$.
Then, an argument on the virtual reduction number allows for the conclusion.

Section~\ref{zeros} treats the case of a sparsing degeneration, where one replaces the entries by zeros in a region in the form of an equilateral triangle all the way to the lower right corner of the matrix. Letting $r$ denote the number of zeros along one side of the triangle, we introduce the number $\mathfrak{o}(r)$ of distinct entries of the generic symmetric matrix of order $m\times m$ set to zero. As it turns, this number will come out as an interesting invariant.
This section is also divided in a similar vein as the previous section, except that the outcome is entirely distinct.
Our first task is to show that the codimension of the gradient ideal is $\leq 3$ (easy) and equals $3$ if and only if $r\leq m-3$ (``if'' is harder).
Next, the Hessian vanishes, so the second matter is to prove that the polar map is a birational map onto the image and to give the nature of the latter.
As it turns, it is a symmetric ladder determinantal ring of dimension ${{m+1}\choose 2}-2\mathfrak o (r)$.
As a consequence, in the second subsection we derive that the image of the (birational) map defined by the submaximal minors is a cone over the polar image with vertex cut out by $\mathfrak o (r)$ coordinate hyperplanes.

Section~\ref{DUAL} deals with the dual variety $V(f)^*$ of $V(f)$.
We first establish a result for a more general sort of degeneration of the generic symmetric $m\times m$ ($m\geq 3$) matrix that includes the two kinds dealt with in the paper. Namely, for such an inclusive setup the dimension of $V(f)^*$ is at least $m-1$.
It will turn up that in both cases of degeneration the dimension is actually $m-1$ -- in line with a point made by Landsberg et. al. in \cite{Landsberg} --
and yet the two cases go apart in their specifics.
Thus, in the cloning case, $V(f)^*$ is properly contained in a ladder determinantal variety of dimension $m-1$, defined by $2\times 2$ minors. The latter is a nice arithmetically Cohen--Macaulay variety, but not so $V(f)^*$ itself whose homogeneous defining ideal contains quadric trinomials as minimal generators.
A supplementary result is that $f$ is a factor of its Hessian determinant with multiplicity ${{m}\choose 2}-2$, equal to the expected multiplicity (\`a la Segre).

As for the zero sparse case, $V(f)^*$ is a ladder determinantal variety of dimension $m-1$, defined by $2\times 2$ minors; in particular, it is an arithmetically Cohen--Macaulay variety and, in addition, it is arithmetically Gorenstein if and only if $r=m-2$.

For the sake of quick browsing, the main results are Theorem~\ref{cloning_sym}, Theorem~\ref{structure_of_submaximal}, Theorem~\ref{polar_structure_zeros}, Theorem\ref{submaximal_are_birational}, Theorem~\ref{dual_ladderlike}, Theorem~\ref{dual_cloning}, Proposition~\ref{expected_multiplicity} and Theorem~\ref{dual_zeros}.

\smallskip

We assume throughout that the ground field has characteristic zero.

\section{Review of basic invariants}\label{basic}

%Let $(R,\mathfrak{m})$ denote a Notherian local ring and its maximal ideal (respectively, a standard graded ring over a field and its irrelevant ideal).
%For an ideal $I\subset \mathfrak{m}$ (respectively, a homogeneous ideal $I\subset \mathfrak{m}$), the \emph{special fiber} of $I$ is the ring $\mathcal{R}(I)/\mathfrak{m}\mathcal{R}(I)$.
%Note that this is an algebra over the residue field of $R$.
%The (Krull) dimension of this algebra is called the \emph{analytic spread} of $I$ and is denoted $\ell (I)$. 

Quite generally, given ideals $J\subset I$  in a ring $R$,  $J$ is said to be a \emph{reduction} of $I$ if there exists an integer $n\geq 0$ such that $I^{n+1}=JI^n.$
An ideal and any of its reductions share the same radical, hence they share the same set of minimal primes and have the same codimension.
A reduction $J$ of $I$ is called \emph{minimal} if no ideal strictly contained in $J$ is a reduction of $I$.
The \emph{reduction number} of $I$ with respect to a reduction $J$ is the minimum integer $n$ such that $JI^{n}=I^{n+1}$. It is denoted by $\mathrm{red}_{J}(I)$. The (absolute) \emph{reduction number} of $I$ is defined as $\mathrm{red}(I)=\mathrm{min}\{\mathrm{red}_{J}(I)~|~J\subset I~\mathrm{is}~\mathrm{a}~\mathrm{minimal}~\mathrm{reduction}~\mathrm{of}~I\}.$

Suppose now that  $R=k[x_0,\ldots,x_d]$ is a standard graded over a field $k$ and $I$ is minimally generated by $n+1$ forms of same degree $s$. In this case, $I$ is more precisely  given by means of a free graded presentation $$ R(-(s+1))^\ell \; \oplus\; \sum_{j\geq 2} R(-(s+j))\stackrel{\varphi}{\longrightarrow} R(-s)^{n+1}\longrightarrow I\longrightarrow 0  $$
for suitable shifts.   Of much interest in this work is the value of $\ell$.  The image of $R(-(s+1))^\ell$ by $\varphi$ is the {\it  linear part of $\varphi$} -- often denoted $\varphi_1$.  It is easy to see that the rank of $\varphi_1$  does not depend  on the particular minimal system of generators of $I$. Thus, we call it the {\it linear rank of} $I$. One says that $I$ has {\it maximal liner  rank} provided its linear rank is $n$ (=rank($\varphi$)). Clearly, the latter condition  is trivially satisfied if $\varphi=\varphi_1 $, in which case $I$ is said  to have {\it linear presentation} (or is {\it linearly presented}). 

Note that $\varphi$  is a graded matrix whose columns generate the (first) {\it syzygy module of $I$} (corresponding to the given choice of generators)  and  a {\it syzyzy} of $I$  is an element of this module -- that is, a linear  relation, with coefficients in $R$, on the chosen generators. In this context, $\varphi_1$ can be taken as the submatrix of $\varphi$  whose entries are linear forms of the standard graded ring $R$.  Thus, the linear rank is the rank of the matrix of the linear syzygies.

A set of $m+1$ forms $f_0,\ldots,f_m$ of the same degree in $R$ defines a rational map $\pp^d\dashrightarrow \pp^m$.
The homogeneous $k$-subalgebra $k[f_0,\ldots,f_m]\subset R$ is up to renormalization isomorphic to the homogeneous defining ring of the image of the map in $\pp^m$.

Given a homogeneous polynomial $f\in R$, its {\em polar map} is the rational map defined by its partial derivatives.
The image of the polar map of $f$ is called its {\em polar variety}.
 $f$ is said to be {\em homaloidal} if its polar map is birational. 
 The ideal of $R$ generated by the partial derivatives of $f$ will often be called the {\em gradient ideal} of $f$.
 
 We refer to \cite{B-H} for other ideal theoretic notions not reviewed in this section and to \cite{B-V} for related ideas on determinantal rings.
 As a guide for some of the algebraic side of birational maps one can look at \cite{AHA}, while the basic material on initial ideals can be traced to \cite{HeHiBook}.

\section{Cloning}\label{CLONING}

A simple coordinate-like degeneration of (\ref{symgeneric}) is obtained via an endomorphism that fixes all variables appearing as entries, except one, and maps the latter to one of the other variables. 
%This procedure  has been dubbed {\em cloning} in \cite{LiZaSi}.

In this section, we deal with the case where the cloning affects two entries on the main diagonal -- we may in this case refer to an $\mathcal{MD}$-{\em cloning} ($\mathcal{MD}$ for ``main diagonal'').
Cloning along anti-diagonals will be considered elsewhere.
Arbitrary cloning may have an unexpected behavior (see Example~\ref{bad_cloning} below) as compared to the $\mathcal{MD}$-cloning.

Clearly, in the case of an $\mathcal{MD}$-cloning, by suitable permutation of rows and columns, without disrupting symmetry, we can move the two entries affected in this process all the way down to the bottom right of the main diagonal. Thus, starting out from $\mathcal S$ as in (\ref{symgeneric}), we may assume that the entry $x_{m,m}$ is replaced  by $x_{m-1,m-1}$, so that the cloned matrix has the form: 

\begin{equation}\label{symgeneric_cloned}
\mathcal{SC}:=\left(
\begin{matrix}
x_{1,1}&x_{1,2}& x_{1,3}&\ldots & x_{1,m-2}  & x_{1,m-1} & x_{1,m}\\
x_{1,2}&x_{2,2}&x_{2,3} &\ldots & x_{2,m-2}  & x_{2,m-1} & x_{2,m}\\
x_{1,3}&x_{2,3}&x_{3,3} &\ldots & x_{3,m-2}  & x_{3,m-1} & x_{3,m}\\
\vdots &\vdots &\ldots &\vdots & \vdots & \vdots & \vdots\\
x_{1,m-2} & x_{2,m-2}& x_{3,m-2}&\ldots & x_{m-2,m-2}  & x_{m-2,m-1} & x_{m-2,m}\\
x_{1,m-1} & x_{2,m-1}& x_{3,m-1} &\ldots & x_{m-2,m-1}  & x_{m-1,m-1} & x_{m-1,m}\\
x_{1,m} & x_{2,m}& x_{3,m} &\ldots & x_{m-2,m}  & x_{m-1,m} & x_{m-1,m-1}
\end{matrix}
\right).
\end{equation}
Note that we have traded the general notation $\mathcal{DS}$ in the Introduction for the present one.
Throughout this part we let $R$ denote the polynomial ring generated over the field $k$ by the entries of $\mathcal{SC}$. For a given integer $1\leq t\leq m$, 
 $I_t({\mathcal{SC}})$ denotes the ideal of $R$ generated by the $t$-minors of $\mathcal{SC}$.

The following basic result is customarily quoted as a consequence of \cite{Golberg}. An independent proof appeared in \cite[Proposition 5.3.1]{Maral}. We restate it noting that a certain hypothesis in the latter is unnecessary and give a short proof for the reader convenience.

\begin{Proposition}\label{GolMar}
	Let ${M}$ denote a square matrix over a polynomial ring $R=k[y_1,\ldots,y_n]$  such that  every entry is either 0 or $y_i$ for some $i=1,\ldots,n$. 
	Then, for each $i=1,\ldots, n$, the partial derivative of  $f=\det({M})$ with respect to $y_i$ is the sum
	of the {\rm (}signed{\rm )} cofactors of $y_i$ in all its slots as an entry of $M$.
\end{Proposition}
\demo
More generally, let $N$ denote an $m\times m$ matrix with linear $l_{i,j}$ entries in the polynomial ring $R=k[y_1,\ldots,y_n]$. Let $G=(x_{i,j})$ stand for the generic $m\times m$ matrix and write $f:=\det N, g=\det G$.
The ordinary chain rule yields for $1\leq r\leq n$:
$$\partial f(y_1,\ldots,y_n)/ \partial y_r = \sum_{1\leq i,j\leq m} (\partial l_{i,j} / \partial y_r)( \partial g/ \partial x_{i,j} )(l_{i,j}),$$
where $( \partial g/ \partial x_{i,j} )(l_{i,j})$ is to be understood as the result of evaluating the polynomial $\partial g/ \partial x_{i,j}$ by $x_{i,j} \mapsto l_{i,j}$.

Now taking $N=M$, the only non-vanishing terms on the right-hand side of the above expression correspond to the entries $l_{i,j}=y_r$.
On the other hand, $\partial g/ \partial x_{i,j}$ is well-known to be the cofactor of $x_{i,j}$ in the generic matrix $G$. Therefore, when $l_{i,j}=y_r$ the summand
$ (\partial l_{i,j} / \partial y_r)( \partial g/ \partial x_{i,j} )(l_{i,j})=( \partial g/ \partial x_{i,j} )(l_{ij})$ is the cofactor of the entry $y_r$ in slot $(i,j)$ of $M$.
\qed

\subsection{$\mathcal{MD}$-cloning: the polar map and homaloidness}

Throughout we set $f:=\det (\mathcal{SC})$ and let $J=J_f\in R$ denote the gradient ideal of $f$, i.e., the ideal generated by the partial derivatives of $f$ with respect to the variables of $R$, the polynomial ring in the entries of $\mathcal{SC}$ over the ground field $k$. 
As usual, $I_t(\mathcal{SC})$ will denote the ideal generated by the $t$-minors of $\mathcal{SC}$.

Let $f_{i,j}$ denote the  $x_{i,j}$-derivative of $f$ and let $\Delta_{j,i}$ stand for the (signed) cofactor of the $(i,j)$th entry of $\mathcal{SC}$. One knows that, by symmetry, the equality $\Delta_{i,j}=\Delta_{j,i}$ holds over the generic symmetric matrix. Since the passage to the cloned version is via a ring homomorphism and the latter commutes with formation of minors, the equality holds over $\mathcal{SC}$ as well. This remark, and possibly others in the same vein, will be used without further ado.

\begin{Theorem}\label{cloning_sym}
Consider the cloned matrix as in {\rm (\ref{symgeneric_cloned})}, with $m\geq 3$.
One has:
\begin{itemize}
\item[{\rm (i)}]   $J$ is a codimension $3$ ideal contained in $I_{m-1}(\mathcal{SC})$.
\item[{\rm (ii)}]  $R/(f)$ is a normal domain.
\item[{\rm (iii)}] The Hessian determinant $h(f)$ does not vanish.
\item[{\rm (iv)}] The linear rank of the gradient ideal of $f$ is ${{m+1}\choose {2}}-2$ {\rm (}maximum possible{\rm )}.
\item[{\rm (v)}] $f$ is homaloidal.
\end{itemize}
\end{Theorem}
\demo
%Let $f_{i,j}$ denote the  $x_{i,j}$-derivative of $f$ and let $\Delta_{j,i}$ stand for the (signed) cofactor of the $(i,j)$th entry of $\mathcal{SC}$. Observe that, by symmetry, one has $\Delta_{i,j}=\Delta_{j,i}$.
(i) By Proposition~\ref{GolMar}, the partial derivative $f_{i,i}$ coincides with the cofactor of $x_{i,i}$, for $i=1,\ldots, m-2$, while $f_{m-1,m-1}$ is the sum of the respective cofactors of $x_{m-1,m-1}$ corresponding to its two slots. By the same token, $\Delta_{j,i}=1/2f_{i,j}$, for all $i\neq j$.
In particular, $J\subset I_{m-1}(\mathcal{SC})$. On the other hand, the ideal of the submaximal minors of the  $m\times m$ generic symmetric matrix specializes since it is a prime ideal generated in degree $\geq 2$. Therefore,  the codimension of $I_{m-1}(\mathcal{SC})$ is $3$ and hence, the codimension of $J$ is at most $3$.

To show that the codimension of $J$ is exactly $3$ we consider the initial ideal of $J$ in the reverse lexicographic order. For $m\geq 5$, direct inspection shows that for $m$ odd one has
\begin{eqnarray*}
	{\rm in}(f_{1,1})&=& x_{2,m}^2\cdot  x_{3,m-1}^2  \cdots  x_{\lfloor\frac{m+2}{2}\rfloor,\lfloor\frac{m+2}{2}\rfloor+1}^2\\
	{\rm in}(f_{1,m})&=&2x_{1,m}\cdot  x_{2,m-1}^2 \cdots  x_{\frac{m+1}{2}-1,\frac{m+1}{2}+1}^2\cdot x_{\frac{m+1}{2},\frac{m+1}{2}}\\
	{\rm in}(f_{m-1,m-1})&=&x_{1,m-1}^2\cdot  x_{2,m-2}^2 \cdots  x_{\lfloor\frac{m}{2}\rfloor,\lfloor\frac{m}{2}\rfloor+1}^2
\end{eqnarray*}
while for $m$ even it obtains
\begin{eqnarray*}
	{\rm in}(f_{1,1})&=&x_{2,m}^2\cdot  x_{3,m-1}^2\cdots x_{\frac{m+2}{2}-1,\frac{m+2}{2}+1}^2\cdot x_{\frac{m+2}{2},\frac{m+2}{2}}\\
	{\rm in}(f_{1,m})&=&2x_{1,m}\cdot  x_{2,m-1}^2\cdots x_{\lfloor\frac{m+1}{2}\rfloor,\lfloor\frac{m+1}{2}\rfloor+1}^2\\
	{\rm in}(f_{m-1,m-1})&=&x_{1,m-1}^2\cdot  x_{2,m-2}^2\cdots x_{\frac{m}{2}-1,\frac{m}{2}+1}^2\cdot x_{\frac{m}{2},\frac{m}{2}}
\end{eqnarray*}
Since there are no common variables among the three monomials in each bloc, it follows that ${\rm in}(J)$ has codimension at least 3.

For $m=3$, an easy verification shows that the monomials    $x_{1,2}^2,x_{2,2}^2$ and $x_{1,3}^3$ belong to ${\rm in}(J)$. For $m=4$, which is the hardest case, we resort to a calculation with \cite{M2} to find a minimal set of generators of ${\rm in}(J)$ as follows:
\begin{eqnarray*}\nonumber
&&	x_{2,3}^2x_{3,3},\; x_{1,4}x_{2,3}x_{2,4},\; x_{1,4}x_{2,3}^2,\; x_{1,4}^2x_{2,3},\; x_{1,4}^2x_{2,2}x_{3,4},\; x_{1,4}^2x_{2,2}x_{3,3},\; x_{1,4}^2x_{2,2}x_{2,4}\\  \nonumber
	&& x_{1,4}^3x_{2,2},\; x_{1,3}x_{2,3}x_{3,3},\; x_{1,3}x_{1,4}x_{2,3},\; x_{1,3}x_{1,4}x_{2,2},\; x_{1,3}^2x_{3,3},\; x_{1,3}^2x_{2,2}.
\end{eqnarray*}
It suffices to observe that no two variables divide simultaneously these monomials.
Alternatively, one can verify that the transposed log matrix of these monomials 

{\small $$\left(
	\begin{matrix}
	0& 0 & 0 & 0 & 0& 2& 0& 1& 0\\
	0& 0 & 0 & 1 & 0& 1& 1& 0& 0\\
	0& 0 & 0 & 1 & 0& 2& 0& 0& 0\\
	0& 0 & 0 & 2 & 0& 1& 0& 0& 0\\
	0& 0 & 0 & 2 & 1& 0& 0& 0& 1\\
	0& 0 & 0 & 2 & 1& 0& 0& 1& 0\\
	0& 0 & 0 & 2 & 1& 0& 1& 0& 0 \\
	0& 0 & 0 & 3 & 1& 0& 0& 0& 0\\
	0& 0 & 1 & 0 & 0& 1& 0& 1& 0\\
	0& 0 & 1 & 1 & 0& 1& 0& 0& 0\\
	0& 0 & 1 & 1 & 1& 0& 0& 0& 0\\
	0& 0 & 2 & 0 & 0& 0& 0& 1& 0\\
	0& 0 & 2 & 0 & 1& 0& 0& 0& 0\\
	\end{matrix}
	\right)$$}
is such that any two columns has a null row.

\medskip

(ii) Since $\codim J=3$ by (i), then $R/(f)$ satisfies the property ($R_1$) of Serre and hence it is normal.  Since $f$ is homogeneous, $R/(f)$ is a domain.

\medskip

(iii) Consider the ring endomorphism $\phi_{\mathbf{v}}$ of $R$ by mapping any variable in  $\mathbf{v}$ to itself and by mapping any variable off  $\mathbf{v}$ to zero, where  $\mathbf{v}:=\{x_{1,1},x_{2,2},x_{3,3},\ldots,x_{m-1,m-1}\}$ is the set of variables along the main diagonal.
Let $\mathcal{H}'$ denote the matrix which results by applying $\phi_{\mathbf{v}}$ to the entries of the Hessian matrix $\mathcal{H}(f)$ of $f$.
Clearly, it suffices to show that $\det \mathcal{H}'\neq 0$.

As already observed, the partial derivative $f_{i,i}$ coincides with the cofactor of $x_{i,i}$, for $i=1,\ldots, m-2$, while $f_{m-1,m-1}$ is the sum of the respective cofactors of $x_{m-1,m-1}$ corresponding to its two slots. 
By expanding each such a cofactor according to the Leibniz rule it is  clear that it has a unique (nonzero) term whose support lies in $\mathbf{v}$ and, moreover,  the remaining terms have degree at least 2 in the variables off $\mathbf{v}$. Slight inspection reveals that in the two cofactors of $x_{m-1,m-1}$ the terms  supported in the variables of $\mathbf{v}$ coincide. 

As for $x_{i,j}\notin \mathbf{v}$,  $f_{i,j}=2\Delta_{i,j}$. The Leibniz expansion of this cofactor has no term with support  in $\mathbf{v}$ and has exactly one nonzero term of degree $1$ in the variables off $\mathbf{v}$.

From these observations follows that applying $\phi_{\mathbf{v}}$ to any second partial derivative of $f$ will return zero or a monomial supported on the variables in $\mathbf{v}$.
Thus, the entries of $\mathcal{H}'$ are either zeros or monomials supported on the variables in $\mathbf{v}$.
	
To see that the determinant of the matrix $\mathcal{H}'$ is nonzero, consider the Jacobian matrix of the set of partial derivatives $\{f_v\,|\, v\in\mathbf{v}\}$ with respect to the variables in $\mathbf{v}$.
Let $M_0$ denote the matrix resulting of applying $\phi_{\mathbf{v}}$ to the entries of this Jacobian matrix, considered as a corresponding submatrix of $\mathcal{H}'$.
Up to permutation of rows and columns of $\mathcal{H}'$, we may write
	$$\mathcal{H}'=
	\left(
	\begin{array}{cc}
	M_0 & N \\
	P & M_1
	\end{array}
	\right),
	$$
	where $M_1$ has exactly one nonzero entry on each row and each column. 
	Now, by the way the second partial derivatives of $f$ map via $\phi_{\mathbf{v}}$, as explained above, one must have $N=P=0$.
	Therefore, $\det(\mathcal{H}')=\det(M_0)\det(M_1)$. 
	It now suffices to verify the nonvanishing of these two subdeterminants.
	
This is clear for $M_1$, since it has exactly one nonzero entry on each row and each column.
As for $M_0$, we see that it is the Hessian matrix of the form
$$g:=\left(\prod_{i=1}^{m-2} x_{i,i}\right) x_{m-1,m-1}^2.
$$

This is the product of the generators of the $k$-subalgebra
$$k[x_{1,1},\ldots,x_{m-2,m-2}, x_{m-1,m-1}^2]\subset k[x_{1,1},\ldots,x_{m-2,m-2}, x_{m-1,m-1}].
$$
Clearly these generators are algebraically independent over $k$, hence the subalgebra is isomorphic to a polynomial ring itself.
%Then $g$ becomes the product of the variables of a polynomial ring over $k$.
This is a classical homaloidal polynomial, hence we are done here too.
	
	\medskip
	
(iv)   Using  the Cauchy cofactor identity
 \begin{equation}\label{CofactorFormula}
\mathcal{SC} \cdot {\rm adj}(\mathcal{SC})={\rm adj}(\mathcal{SC})\cdot\mathcal{SC}=\det(\mathcal{SC})\, \mathbb{I}_m
\end{equation}
we find  the following linear relations involving the cofactors  of  $\mathcal{SC}$:

	\begin{equation}\label{sc1} 
	\sum_{j=1}^m \widehat{x_{i,j}}\Delta_{j,1}=0,\ \mbox{for} \; 2\leq i \leq m-1;
	\end{equation}
	\begin{equation}\label{sc2}
\sum_{j=1}^m \widehat{x_{i,j}}\Delta_{j,k}=0,\ \mbox{for} \; 2\leq k  \leq m-2\; \mbox{and} \; k-1\leq i \leq m-1 \; ( k\neq i);
	\end{equation}
	\begin{equation}\label{sc3}
	\sum_{j=1}^{m-1} \widehat{x_{m,j}}\Delta_{j,k}+x_{m-1,m-1}\Delta_{m,k}=0,\ \mbox{for} \; 1\leq k  \leq m-2;
	\end{equation}
	
	\begin{equation}\label{sc4}
	\sum_{i=1}^{m-1} x_{i,m-1} \Delta_{i,m} +x_{m-1,m} \Delta_{m,m}=0
	\end{equation}

	\begin{equation}\label{sc5}
	\sum_{i=1}^{m-2} x_{i,m}\Delta_{i,m-1}+x_{m-1,m}\Delta_{m-1,m-1}+x_{m-1,m-1} \Delta_{m,m-1}=0\end{equation}
	\begin{equation}\label{sc6}
	\sum_{i=1}^{m-2} x_{i,m-1}\Delta_{i,m-1}+x_{m-1,m-1}\Delta_{m-1,m-1} + x_{m-1,m}\Delta_{m,m-1}=\sum_{j=1}^m x_{1, j}\Delta_{j,1}\end{equation}
	\begin{equation}\label{sc7}\sum_{j=1}^{m-1} x_{j,m}\Delta_{j,m}+x_{m-1,m-1}\Delta_{m,m}=\sum_{j=1}^m x_{1, j}\Delta_{j,1}.
	\end{equation}
Here we have set 
\begin{equation}\label{defining_hat}
\widehat{x_{i,j}}=\left\{\begin{array}{ll}
x_{i,j}\quad  \text{if $i\leq j$}\\
x_{j,i} \quad \text{if $i\geq j$}.
\end{array}
\right.
\end{equation}

This notation will be used throughout on several occasions.

Since, as already remarked, one has $f_{i,j}=2\Delta_{i,j}$ for $1\leq i<j\leq m$ and $f_{i,i}=\Delta_{i,i}$ for $1\leq i\leq m-2$, then $(\ref{sc1})$, $(\ref{sc2})$ and $(\ref{sc3})$ give linear syzygies of the partial derivatives. Moreover, since  $f_{m-1, m-1}=\Delta_{m-1,m-1}+\Delta_{m,m}$, adding (\ref{sc4}) to (\ref{sc5}) and (\ref{sc6}) to (\ref{sc7}) outputs two additional linear syzygies of the partial derivatives of $f$.
Thus one has counted a total of $(m-1)+(m-1)+(m-2)+\ldots +3+2={{m+1}\choose {2}}-2$ linear syzygies of $J$.
In order to see that they are moreover independent, we order the set of partial derivatives $f_{i,j}$ in accordance with the following ordered list of the entries $x_{i,j}$:
	\begin{eqnarray} \nonumber
	x_{1,1}, x_{1,2},\ldots, x_{1,m}\rightsquigarrow x_{2,2},\ldots, x_{2,m}\rightsquigarrow\ldots\rightsquigarrow x_{m-2,m-2},x_{m-2,m-1} x_{m-2,m} \rightsquigarrow  x_{m-1,m-1},x_{m-1,m}.
	\end{eqnarray}
	
We now claim that, ordering the set of partial derivatives $f_{i,j}$ in this way, the above sets of linear relations can be grouped into the following block matrix of  linear syzygies:

	$$\left(\begin{array}{cccccc|ccc}
	\Phi_1 &   \ldots &         &                        &   &                                 &\\
	0^{m-1}_{m-1}      &\Phi_2 &\ldots   &      &                    &                              &\\
	0^{m-1}_{m-2}      &0^{m-1}_{m-2} & \Phi_3   &      &                    &                              &\\
	
	\vdots    &\vdots    & \vdots  & \ddots  &  &              &                     &\\
	0 ^{m-1}_4      & 0^{m-1}_4        & 0^{m-2}_4 & \ldots     & \Phi_{m-3}   &       &                        &\\
	0 ^{m-1}_3       & 0^{m-1}_3        &  0^{m-2}_3& \ldots    & 0_3^4&  \Phi_{m-2}          &                        &\\
	\hline
	0^{m-1}_1   & 0_1^{m-1}        & 0_1^{m-2} & \ldots  &   0_1^{4}  &     0_1^{3}       &   x_{m-1,m}& x_{m-1,m-1} \\
	0^{m-1}_1    & 0_1^{m-1}       & 0_1^{m-2}  &  \ldots   &0_1^{4}   & 0_1^3        &   x_{m-1,m-1} & x_{m-1,m}\\
	\end{array}
	\right),$$
	
	where:
	\begin{itemize}
		\item $\Phi_1$ is the matrix obtained  from $(\mathcal{SC})^t$ obtained by multiplying the first row by 2 and omitting the first column.
		\item $\Phi_2$ is the matrix obtained  from $(\mathcal{SC})^t$ by  multiplying the second row by $2$ and omitting the second column and the first row .
		\item  For  $l=3,\ldots,m-2$, $\Phi_l$ is  the matrix obtained  from $(\mathcal{SC})^t$ by multiplying the $l$th row by $2$ and omitting the columns $1,\ldots,l-2,l$ and the rows $1,\ldots,l-1$. 
		\item $0_r^c$ denotes  a zero block of size $r\times c $.
	\end{itemize}
	
Justification is as follows.
	
First, as already observed,  the relations (\ref{sc1}) through (\ref{sc7}) yield linear syzygies of the partial derivatives of $f$.
	
Using the relation between partial derivatives and cofactors,  $(\ref{sc1})$ can be written as 
$$2\,\widehat{x_{i,1}}\,f_{1,1} +\sum_{j=2}^m \widehat{x_{i,j}}\,f_{1,j}=0 ,$$ for $i=2,\ldots , m-1$. Moreover, setting $k=1$ in $(\ref{sc3})$ yields
$$ 2\,\widehat{x_{m,1}}\,f_{1,1}+\sum_{j=2}^{m-1} \widehat{x_{m,j}}\,f_{1,j}+x_{m-1,m-1}\,f_{1,m}=0$$ 
	
Ordering the set of partial derivatives $f_{i,j}$ as explained before, and trading the coefficients of these relations back to variable notation, one gets 
	
	\begin{equation*}\nonumber
	\Phi_1:=\left(\begin{array}{ccccc}
		2x_{1,2} & 2x_{1,3} & \ldots  &2x_{1,m-1}  & 2x_{1,m} \\
		x_{2,2} &x_{2,3} & \ldots  &x_{2,m-1}  & x_{2,m} \\
		\vdots  & \vdots & \cdots  &\vdots  & \vdots \\
		x_{2,m-1} & x_{3,m-1}& \ldots  &x_{m-1, m-1} & x_{m-1,m}\\
		x_{2,m} & x_{3,m}& \ldots  &x_{m-1, m} & x_{m-1,m-1}
		\end{array}
		\right).
			\end{equation*}

Note that $\Phi_1$ coincides indeed with the submatrix of $\mathcal{SC}^t$ obtained by multiplying the first row by 2 and omitting the first column.
	
Continuing, for each $l=2,\ldots ,m-2$  the  block $\Phi_l$ comes from the relation (\ref{sc2})   and (\ref{sc3}) (setting $k=l$).  Finally, the lower right corner $2\times 2$ block of the matrix of  linear syzygies comes from the last two relations obtained by adding (\ref{sc4}) to (\ref{sc5}) and  (\ref{sc6}) to (\ref{sc7}).
	
So much for justification.

Now, counting through the sizes of the various blocks, one sees that this matrix is $({{m+1}\choose {2}}-1)\times ({{m+1}\choose {2}}-2)$.
Omitting its first row obtains a block-diagonal submatrix of size
$({{m+1}\choose {2}}-2)\times ({{m+1}\choose {2}}-2)$, where each block has nonzero determinant. Thus, the linear rank of $J$ attains the maximum.

	\medskip

(v) By (iii) the polar map of $f$ is dominant. Since the linear rank is maximum by (iv), one can apply  \cite[Theorem 3.2]{AHA} to conclude that $f$ is homaloidal.
\qed

\medskip

The following example shows that an arbitrary cloning may lack most of the properties listed in Theorem~\ref{cloning_sym}.

\begin{Example}\label{bad_cloning}\rm
Consider the cloning endomorphism on the $3\times 3$ generic symmetric matrix that maps $x_{2,3}$ to $x_{1,1}$. Changing the names of the variables for the sake of visualization, the resulting degeneration is the matrix
$$\mathfrak{S}=
\left(
\begin{matrix}
x_1 & x_2 & x_3\\
x_2 & x_4 & x_1\\
x_3 & x_1 & x_5
\end{matrix}
\right).
$$
A calculation with \cite{M2} gives that the linear rank of the gradient ideal $J\subset k[x_1,x_2,x_3,x_4,x_5]$ of $\det \mathfrak{S}$ is $3$,  one below the maximum $5-1=4$. Thus, item (iv) of Theorem~\ref{cloning_sym} fails here. In fact, another computation with \cite{M2} yields that $J$ is an ideal of linear type, hence by \cite[Proposition 3.4]{AHA}, $\det \mathfrak{S}$ is not homaloidal, showing that item Theorem~\ref{cloning_sym} (v) fails as well.
Note that this degeneration has a homological behavior reminiscent of the $3\times 3$ generic Hankel matrix (see \cite[Proposition 3.21 (b)]{MAron}).
\end{Example}

\subsection{The structure of the submaximal minors}

In this part we study the nature of the ideal of submaximal minors of $\mathcal{SC}$. As previously, $J$ denotes the gradient ideal of $f=\det \mathcal{SC}$.

\begin{Theorem}\label{structure_of_submaximal}
	Consider the matrix $\mathcal{SC}$ as in $(\ref{symgeneric_cloned})$, with $m\geq 4$.  Set $P:=I_{m-1}(\mathcal{SC})\subset R$.
	Then 
	\begin{enumerate}
		\item[{\rm (i)}] $R/P$ is a Cohen-Macaulay normal domain of codimension $3$.
		\item[{\rm (ii)}] $P$ is the minimal component of the primary decomposition of $J$ in $R$.
		\item[\rm (iii)] $J$ defines a double structure on the variety defined by $P$, with a unique embedded component of codimension $2(m-1)$ supported on a linear space and no other embedded component of codimension $\leq 2(m-1)$.
		\item[\rm (iv)] The $(m-1)$-minors of $\mathcal{SC}$ define a birational map $\pp^{{{m+1}\choose {2}}-2}\dasharrow \pp^{{{m+1}\choose {2}}-1}$ onto a hypersurface of degree $m-1$ with defining equation $\mathbb{D}_{m,m}-\mathbb{D}_{m-1,m-1}$, where $\mathbb{D}_{m,m}$ and $\mathbb{D}_{m-1,m-1}$  denote  the cofactors of $y_{m,m}$  and $y_{m-1,m-1}$, respectively, in the $m\times m$ generic symmetric 
		matrix $(y_{i,j})_{1\leq i\leq j\leq m}$ on the target coordinates.	
		\item[\rm (v)] 	$J$ is not a reduction of $P$.
	\end{enumerate}
\end{Theorem}
\demo 
(i) Since $P$ is a specialization of the corresponding ideal of the generic symmetric matrix, it follows that $R/P$ is a Cohen-Macaulay ring of codimension $3$. 
 
As $P$ is homogeneous, normality of $R/P$ implies that $P$ is prime. To show normality,  Serre's property $(S_2)$ is automatic since $R/P$ is Cohen-Macaulay. Therefore,  it remains to prove that it satisfies condition  $(R_1)$. For this, let $\Theta$ denote the Jacobian matrix of the generators of $P$ with respect to the variables of $R$. We proceed to show that  
$$\codim(I_3 (\Theta),P) \geq 5.$$ 

We will argue via the initial ideal in the revlex monomial order induced by the ordering the variables in the sequence in which they appear in the matrix respecting the rows. 

For $m = 4$, direct inspection shows that the monomials $x_{3,3}^6$, $ x_{2,4}^6$, $x_{2,3}^6$, $x_{1,4}^6$ and $x_{1,3}^6$ belong to ${\rm in}(I_3(\Theta))$.  For $m = 5$, inspection is harder so we resort to a calculation with \cite{M2} to find the following  monomials in the initial ideal of $(I_3 (\Theta),P)$: $x_{1,3}^6x_{2,2}^3x_{4,5}$, $x_{1,4}^2x_{2,3}^2$, $x_{1,5}^3x_{2,4}^6$, $x_{2,5}^6x_{3,3}^3$ and $x_{3,4}^6x_{4,4}^3$. So, for $m=4,5$ we have $\codim(I_3 (\Theta),P) \geq 5$ and thereby $R/P$ is normal.  	

Now assume that $m\geq 6$.

In the generic symmetric case, the leading term of any minor determinant is well-known to be the product of the entries along the its main anti-diagonal. We claim that this remains true for the $(m-1)$-minors and the $(m-2)$-minors of $\mathcal{SC}$, provided $m\geq 6$.  Indeed,  let $M$ denote any such minor and let $D$ denote the  product of the entries along the its main anti-diagonal. Observe that each variable $x_{i,j}$ in the anti-diagonal of $M$ satisfies $x_{i,j}\geq x_{(m+2)/2, (m+2)/2}$ if $m$ is even  (respectively, $x_{i,j}\geq x_{(m+3)/2, (m+3)/2}$ if $m$ is odd). Thus, when $m\geq 6$ one has $x_{i,j}>x_{m-1,m-1}$ since $x_{(m+2)/2, (m+2)/2}>x_{m-1,m-1}$ if $m$ is even (respectively, $ x_{(m+3)/2, (m+3)/2}>x_{m-1,m-1}$ if $m$ is odd). This ensures that, in the revlex monomial order, any monomial involving the cloned variable $x_{m-1,m-1}$ is smaller than $D$ and, therefore, $D$ is the leading term of $\det M$.

To show that $\codim(I_3 (\Theta) ,P) \geq 5$  consider the following submatrices of $\Theta$:

$$ \Theta_1:=\left(\begin{array}{ccc}
\partial \Delta_{2,2}/\partial x_{1,1}  &\partial \Delta_{1,2}/\partial x_{1,1}  &\partial \Delta_{1,1}/\partial x_{1,1}  \\
\partial \Delta_{2,2}/\partial x_{1,2}  &\partial \Delta_{1,2}/\partial x_{1,2}  &\partial \Delta_{1,1}/\partial x_{1,2}  \\
\partial \Delta_{2,2}/\partial x_{2,2}  &\partial \Delta_{1,2}/\partial x_{2,2}  &\partial \Delta_{1,1}/\partial x_{2,2}  \\
\end{array}
\right)$$ 
$$ \Theta_2:=\left(\begin{array}{ccc}
\partial \Delta_{m,m}/\partial x_{1,1}  &\partial \Delta_{1,m}/\partial x_{1,1}  &\partial \Delta_{1,1}/\partial x_{1,1}  \\
\partial \Delta_{m,m}/\partial x_{1,m}  &\partial \Delta_{1,m}/\partial x_{1,m}  &\partial \Delta_{1,1}/\partial x_{1,m} \\
\partial \Delta_{m,m}/\partial x_{m-1,m-1}  &\partial \Delta_{1,m}/\partial x_{m-1,m-1}  &\partial \Delta_{1,1}/\partial x_{m-1,m-1} 
\end{array}
\right)$$ 
and 
$$\Theta_3:= \left(\begin{array}{ccc}
\partial \Delta_{m-1,m}/\partial x_{m-2,m}  &\partial \Delta_{m-1,m-1}/\partial x_{m-2,m}  &\partial \Delta_{m-2,m}/\partial x_{m-2,m}  \\
\partial \Delta_{m-1,m}/\partial x_{m-1,m-1}  &\partial \Delta_{m-1,m-1}/\partial x_{m-1,m-1}  &\partial \Delta_{m-2,m}/\partial x_{m-1,m-1} \\
\partial \Delta_{m-1,m}/\partial x_{m-1,m}  &\partial \Delta_{m-1,m-1}/\partial x_{m-1,m}  &\partial \Delta_{m-2,m}/\partial x_{m-1,m} 
\end{array}
\right)$$ 
The objective is to write the determinants of these submatrices in terms of the determinants of certain $(m-2)\times (m-2)$ submatrices of $\mathcal{SC}$ itself.
The ones we need are as follows:
\begin{equation*}
M_1= \left(  \begin{matrix}  
	x_{3,3} &\ldots & x_{3,m-2}  & x_{3,m-1} & x_{3,m}\\
	\ldots &\vdots & \vdots & \vdots & \vdots\\
	x_{3,m-2}&\ldots & x_{m-2,m-2}  & x_{m-2,m-1} & x_{m-2,m}\\
	x_{3,m-1} &\ldots & x_{m-2,m-1}  & x_{m-1,m-1} & x_{m-1,m}\\
	x_{3,m} &\ldots & x_{m-2,m}  & x_{m-1,m} & x_{m-1,m-1}
\end{matrix}
\right),
\end{equation*}
omitting the first two rows and first two columns.
\begin{equation*}
M_2= \left(  \begin{matrix}  	
	x_{2,2} &\ldots & x_{2,m-2}  & x_{2,m-1} \\
	\vdots &\vdots & \vdots & \vdots \\
	x_{2,m-2}&\ldots & x_{m-2,m-2}  & x_{m-2,m-1} \\
	x_{2,m-1} &\ldots & x_{m-2,m-1}  & x_{m-1,m-1} 
\end{matrix}
\right),
\end{equation*}
omitting the first row and column and the $m$th row and column.
\begin{equation*}
 M_3= \left(  \begin{matrix}  	
	x_{2,2} &\ldots & x_{2,m-2}  & x_{2,m} \\
	\vdots &\vdots & \vdots & \vdots \\
	x_{2,m-2}&\ldots & x_{m-2,m-2}  & x_{m-2,m} \\
	x_{2,m} &\ldots & x_{m-2,m}  & x_{m-1,m-1} 
\end{matrix}
\right)
\end{equation*}
omitting the first row and column and the $(m-1)$th row and column.
$$M_4= \left(  \begin{matrix}  
	x_{1,1} &\ldots & x_{1,m-3}  & x_{1,m-2} \\	
	\vdots &\vdots & \vdots & \vdots \\
	x_{1,m-3}&\ldots & x_{m-3,m-3}  & x_{m-3,m-2} \\
	x_{1,m-1} &\ldots & x_{m-3,m-1}  & x_{m-2,m-1} 
\end{matrix}
\right),
$$
omitting the $(m-2)$th and $m$th rows, and the last two columns.
$$ M_5= \left(  \begin{matrix}  	
	x_{1,1} &\ldots & x_{1,m-3}  & x_{1,m-2} \\	
	\vdots &\vdots & \vdots & \vdots \\
	x_{1,m-3}&\ldots & x_{m-3,m-3}  & x_{m-3,m-2} \\
	x_{1,m-2} &\ldots & x_{m-3,m-2}  & x_{m-2,m-2} 
\end{matrix}
\right),
$$
omitting the last two rows and the last two columns.
$$ M_6= \left(  \begin{matrix}  
	x_{1,1} &\ldots & x_{1,m-3}  & x_{1,m-2} \\	
	\vdots &\vdots & \vdots & \vdots \\
	x_{1,m-3}&\ldots & x_{m-3,m-3}  & x_{m-3,m-2} \\
	x_{1,m} &\ldots & x_{m-3,m}  & x_{m-2,m} 
\end{matrix}
\right),
$$
omitting the $(m-2)$th and $(m-1)$th rows, and the last two columns.
$$ M_7= \left(  \begin{matrix}  	
	x_{1,1} &\ldots & x_{1,m-3}  & x_{1,m-1} \\	
	\vdots &\vdots & \vdots & \vdots \\
	x_{1,m-3}&\ldots & x_{m-3,m-3}  & x_{m-3,m-1} \\
	x_{1,m-1} &\ldots & x_{m-3,m-1}  & x_{m-1,m-1} 
\end{matrix}
\right),
$$
omitting the $(m-2)$th and $m$th rows, and the $(m-2)$th and $m$th columns.

\smallskip

We now analyse the partial derivatives of the various cofactors.

By close inspection of the  cofactors  $\Delta_{1,1},\,\Delta_{1,2},\, \Delta_{2,2}, \Delta_{1,m}, \, \Delta_{m,m}$ one has
$$\partial \Delta_{1,2}/\partial x_{1,1}=\partial \Delta_{1,1}/\partial x_{1,1}=\partial \Delta_{1,1}/\partial x_{1,2}=\partial \Delta_{1,m}/\partial x_{1,1}=\partial \Delta_{1,1}/\partial x_{1,m}=0.
$$ 
Similarly, $\partial\Delta_{m-1,m}/\partial x_{m-1,m-1}=\partial \Delta_{m-1,m-1}/\partial x_{m-1,m}=0$.

Given indices $i,j$, we write $M_{i,j}$ for the matrix such that $\Delta_{i,j}=(-1)^{i+j}\det M_{i,j}$. By Proposition~\ref{GolMar}, for a variable $x_{k,l}$ which is an entry of $M_{i,j}$,   $\partial \Delta_{i,j}/\partial x_{k,l}$  is the sum of the {\rm (}signed{\rm )} cofactors on $M_{i,j}$  of the entry $x_{k,l}$, in all its slots in $M_{i,j}$. 
We thus get 
$$\partial \Delta_{2,2}/\partial x_{1,1}=-\partial \Delta_{1,2}/\partial x_{1,2}=\partial \Delta_{1,1}/\partial x_{2,2}=\det M_1$$
$$\partial \Delta_{m,m}/\partial x_{1,1}= -\partial \Delta_{1,m}/\partial x_{1,m} =\det (M_2)$$ 
$$\partial \Delta_{1,1}/\partial x_{m-1,m-1}=\det (M_2)+\det(M_3)$$

By a similar token, we still get  
$$\partial   \Delta_{m-1,m}/\partial x_{m-2,m}=\partial \Delta_{m-2,m}/\partial x_{m-1,m}=\det(M_4) $$
$$-\partial \Delta_{m-1,m}/\partial x_{m-1,m}= \partial \Delta_{m-1,m-1}/\partial x_{m-1,m-1}=\det(M_5)$$ 
$$\partial \Delta_{m-1,m-1}/\partial x_{m-2,m}=2\,\partial \Delta_{m-2,m}/\partial x_{m-1,m-1}=-2\det(M_6);$$
$$\partial \Delta_{m-2,m}/\partial x_{m-2,m}=-\det(M_7).$$

Collecting these data, one finds:
$$ \det (\Theta_1)=\det\left(\begin{array}{ccc}
\det (M_1)  & 0  & 0  \\
\partial \Delta_{2,2}/\partial x_{1,2}  &-\det (M_1)   & 0  \\
\partial \Delta_{2,2}/\partial x_{2,2}  &\partial \Delta_{1,2}/\partial x_{2,2}  &\det (M_1)   \\
\end{array}
\right)=- \det (M_1)^3,
$$   
   
\begin{eqnarray}\nonumber
\det(\Theta_2)&=&\det {\left(\begin{array}{ccc}
	\det (M_2)  & 0  & 0  \\
	\partial \Delta_{m,m}/\partial x_{1,m}  &-\det (M_2)   & 0  \\
	\partial \Delta_{m,m}/\partial x_{m,m}  &\partial \Delta_{1,m}/\partial x_{m,m}  &\det (M_2) + \det(M_3)   \\
	\end{array}
	\right)}\\ \nonumber
&=&-\det (M_2)^3-\det (M_2)^2\det(M_3).
\end{eqnarray} 
   
\begin{eqnarray}\nonumber \det(\Theta_3) &=& \left(\begin{array}{ccc}
\det(M_4)  &-2\det(M_6)  &-\det(M_7)  \\
0 &\det(M_5)  & -\det(M_6) \\
-\det(M_5)  &0  &\det(M_4) 
\end{array}
\right)\\\nonumber
& &=\det(M_4)^2\det(M_5)-\det(M_5)^2\det(M_7)- 2\det(M_6)^2\det(M_5).
\end{eqnarray}

To complete the argument, we show:

{\sc Claim:} The ideal  $\left( \Delta_{1,1},\, \Delta_{m,m},\, \det(\Theta_1),\,\det(\Theta_2),\,\det(\Theta_3) \right) \subset (I_3 (\Theta) , P)$ has codimension $5$.

 For this, we look at its initial ideal in the revlex order.  According to a remark at the beginning of the proof, when $m\geq 6$ the initial term of an $(m-2)$-minor or an $(m-1)$-minor is the product of the entries along  its main anti-diagonal. Thus,  one gets immediately 
\begin{eqnarray*}\nonumber
 {\rm in}  (\Delta_{1,1}) &=&-\prod_{i+j=m+2}\widehat{x_{i,j}}\\ \nonumber
{\rm in}  (\Delta_{m,m}) &=&-\prod_{i+j=m}\widehat{x_{i,j}}\\  \nonumber
{\rm in}  ( \det (\Theta_1)) &=&({\rm in} (\det (M_1)))^3=-\left( \prod_{i+j=m+3}\widehat{x_{i,j}}\right) ^3
\end{eqnarray*}  

Now, one has $ {\rm in}  ( \det (\Theta_2)) =-{\rm in} ( \det (M_2)^2) \cdot\max\left\lbrace {\rm in} (\det (M_2)),{\rm in}(\det(M_3))\right\rbrace.$

Observing that 
$${\rm in} (\det (M_2))= -\prod_{i+j=m+1, i\neq 1,2}\widehat{x_{i,j}}x_{2,m-1}^2, \quad {\rm in} (\det (M_3)) =-\prod_{i+j=m+1, i\neq 1,2}\widehat{x_{i,j}}x_{2,m}^2$$ 
and since in the revlex monomial order  $x_{2,m}$ is smaller than $x_{2,m-1}$, we conclude that  $ {\rm in}  ( \det (\Theta_2))=( \prod_{i+j=m+1, i\neq 1}\widehat{x_{i,j}})^3$.

Finally, we consider $$\det(\Theta_3)=-2\det(M_6)^2\det(M_5)+ \det(M_4)^2\det(M_5)-\det(M_5)^2\det(M_7).$$ 
Let $D$ denote the product of all variables along the main diagonal of $M_5$ excluding the variables $x_{1,m-2}$, that is, $D=\prod_{i+j=m-1, i\neq 1}\widehat{x_{i,j}}$. We observe  that  \begin{itemize}
	\item ${\rm in}(\det(M_4))= - x_{1,m-2}x_{1,m-1}D$;
	\item  $ {\rm in}(\det(M_5))=  -x_{1,m-2}^2D$;
	\item ${\rm in}(\det(M_6))=  -x_{1,m-2}x_{1,m}D$;
	\item ${\rm in}(\det(M_7))=  -x_{1,m-1}^2D$.
\end{itemize}

Thus,    $${\rm in}(-2\det(M_6)^2\det(M_5))= 2x_{1,m-2}^4x_{1,m}^2D^3$$  and $${\rm in}(\det(M_4))^2\det(M_5))=-x_{1,m-2}^4x_{1,m-1}^2D^3  =(\det(M_5))^2\det(M_7). $$ 
Consequently,  ${\rm in}(\det(M_4)^2\det(M_5)-\det(M_5)^2\det(M_7)) < x_{1,m-2}^4x_{1,m-1}^2D^3$ and    
$$ {\rm in}(\det(\Theta_3))\leq \max\left\lbrace 2x_{1,m-2}^4x_{1,m}^2D^3,\, {\rm in}(\det(M_4)^2\det(M_5)-\det(M_5)^2\det(M_7)) \right\rbrace.$$ 

As  the variable $x_{1,m}$ is absent in the submatrices $M_4$, $M_5$ and $M_6$ of $\mathcal{SC}$, then the support terms of the polynomial  $(\det(M_4))^2\det(M_5)-(\det(M_5))^2\det(M_7)$ do not include the monomials  $x_{1,m-2}^4x_{1,m-1}x_{1,m}D^3$ or  $x_{1,m-2}^4x_{1,m}^2D^3$. This ensures that ${\rm in}(\det(\Theta_3)) = 2x_{1,m-2}^4x_{1,m}^2D^3$.

Since there are no common variables among the monomials 
$${\rm in}  (\Delta_{1,1}) =-\prod_{i+j=m+2}x_{i,j}, \quad {\rm in}  (\Delta_{m,m}) =-\prod_{i+j=m}x_{i,j} , \quad
{\rm in}  ( \det (\Theta_1)) =-\left( \prod_{i+j=m+3}x_{i,j}\right) ^3,$$ $$ {\rm in}  ( \det (\Theta_2))=\left( \prod_{i+j=m+1, i\neq 1}x_{i,j}\right)^3 \quad \mbox{and}\quad    {\rm in}(\det(\Theta_3)) = 2x_{1,m-2}^4x_{1,m}^2D^3,$$ 
it follows that $\left( \Delta_{1,1},\, \Delta_{m,m},\, \det(\Theta_1),\,\det(\Theta_2),\,\det(\Theta_3) \right)$ has codimension $5$ as claimed.

Therefore,  $(I_3 (\Theta) , P)$ has codimension at least $5$, as was to be shown.

\medskip

(ii)  By item (i), $P$ is a prime ideal  of codimension 3. We first show that 
$\codim (J:P)>3$, which ensures that the radical of the unmixed part  of $J$ has no primes of codimension $<3$  and coincides with $P$.   

%Here, the partial derivatives, except $f_{m-1,m-1}$   which coincides with the sum of cofactors $\Delta_{m-1,m-1}+\Delta_{m,m}$,  are scalar multiple of the cofactors and so we can write  $P=(J,\Delta_{m,m})$ and $P= (J,\Delta_{m-1,m-1})$ as in the Theorem~\ref{primality_generic_cloned} and both ways of writing  $P$ will be needed to complete $\codim (J:P)>3$. 

Since $f_{m-1,m-1}=\Delta_{m-1,m-1}+\Delta_{m,m}$ and any other partial derivative  $f_{i,j}$ coincides with $\Delta_{i,j}$, up to escalar multiple, we can write 
$P=(J,\Delta_{m,m})$ and $P= (J,\Delta_{m-1,m-1})$. In particular $J:P= J:\Delta_{m,m}$ and $J:P=J:\Delta_{m-1,m-1}$.

The  cofactor identity yields the following relations:
\begin{eqnarray} \nonumber
\sum_{j=1}^m \widehat{x_{k,j}}\Delta_{j,m}&=&0,\; \mbox{for} \;  k=1,\ldots, m-1\;\\\nonumber
\sum_{i=1}^{m-1} x_{i,m}\Delta_{i,m}+x_{m-1,m-1}\Delta_{m,m}&=&\sum_{j=1}^m x_{1,j}\Delta_{1,j}; \\ \nonumber
%\sum_{i=1}^m x_{i,k}\Delta_{i,m}&=&0, \; \mbox{for} \;  k = m-2,m-1 ;
\end{eqnarray}
where, $\widehat{x_{i,j}}$  is as in (\ref{defining_hat}).

%Since  $f_{i,j}=\Delta_{i,j}$ or $1/2 f_{i,j}=\Delta_{i,j}$  for $(i,j)\neq (m-1,m-1)$, 

By the preceding observation, the above relations imply that the entries of the $m$-th column of $\mathcal{SC}$ belong to the ideal $J:\Delta_{m,m} =J:P $.

In addition, from  the  cofactor identity  we read the following relations:

\begin{eqnarray} \nonumber
\sum_{j=1, j\neq m-1}^m \widehat{x_{k,j}}\Delta_{j,m-1}+ \widehat{x_{k,m-1}}\Delta_{m-1,m-1}&=&0,\; \mbox{for} \;  k=1,\ldots, m, \;( k\neq m-1);\\\nonumber
\sum_{j=1, j\neq m-1}^m \widehat{x_{m-1,j}}\Delta_{j,m-1}+x_{m-1,m-1}\Delta_{m-1,m-1}&=&\sum_{j=1}^m x_{1,j}\Delta_{j,1}; \\ \nonumber
% \sum_{i=1,i\neq m-1}^m x_{i,k}\Delta_{m-1,i}+x_{m-1,k}\Delta_{m-1,m-1}&=&0, \; \mbox{for} \;  k = 1,...,m \;(k\neq m-1) ;
\end{eqnarray}

Then by the same token as above, the entries of the $(m-1)$-th column of $\mathcal{SC}$ belong to the ideal $J:\Delta_{m-1,m-1}=J:P$. 

From this, the variables of the two rightmost columns of $\mathcal{SC}$  conduct $P$ into $J$. In, particular, the codimension of $J:P$  is at least 4, as needed. 
 
Now, since $P$  has codimension 3 then $J:P\not\subset P$. Picking a element $a\in J:P\setminus P$ shows that $P_P\subset J_P$. Therefore $P$ is the unmixed part of $J$.

To  conclude  that $P$ is the minimal primary component of $J$,  we observe that,  by symmetry,  the entries of the last two columns are the same as those  last rows. As is clear that $P$ is contained in the ideal generated by these variables it follows that $P^2\subset J$. 
Therefore, the radical of $J$ -- i.e., the radical of the minimal primary part of $J$-- is $P$.

\medskip

(iii)  By (ii), $P$ is the minimal component of a primary decomposition of $J$ and $I\subset J:P$, where $I$  denotes the ideal  generated by the variables  of the two rightmost columns of $\mathcal{SC}$ -- note that there are $m+m-2=2(m-1)$ such variables. 

{\sc Claim:} $J:P=I$.

\smallskip

To see this, introduce the subideals  $I'\subset I$ (respectively, $I''\subset I$)  generated by the variables on the $(m-1)$th column (respectively, by the variables on the $m$th column).
Clearly, $I=I'+I''$.

Now, $\Delta_{m,m}\notin I''$, while $\Delta_{i,j}\in I''$ for all $(i,j)\neq (m,m)$, since $\Delta_{i,j}$ includes a row or a column with entries in $I''$.
Similarly, $\Delta_{m-1,m-1}\notin I'$, while $\Delta_{i,j}\in I'$ for all $(i,j)\neq (m-1,m-1)$.

Now, recall that $J$ is generated by the cofactors 
$$\Delta_{l,h},\; \mbox{\rm with} \;(l,h)\neq (m-1,m-1), (l,h)\neq (m,m)
$$
and the additional form $\Delta_{m,m}+\Delta_{m-1,m-1}$.
Therefore, given, say,  $b\in J:P=J:  \Delta_{m,m}$, one has
\begin{equation} \label{Q1}
b\,\Delta_{m,m}=\sum_{(i,j)\neq(m-1,m-1)} a_{i,j}\Delta_{j,i}+a(\Delta_{m-1,m-1}+\Delta_{m,m})
\end{equation}
for certain $a_{i,j},a\in R$. 
Then 
$$ (b-a)\Delta_{m,m} =\sum_{(i,j)\neq(m-1,m-1)} a_{i,j}\Delta_{j,i}+a\Delta_{m-1,m-1}  \in I''.$$ 
Since $I''$ is a prime ideal and $\Delta_{m,m}\notin I''$, we have $c:=b-a\in I''$.  
Substituting for $a=b-c$ in  (\ref{Q1}) gives   
$$  (-b+c)\Delta_{m-1,m-1}=\sum_{(i,j)\neq(m-1,m-1)} a_{i,j}\Delta_{j,i}-c\Delta_{m,m}\in I'.$$    

By a similar token, since $\Delta_{m-1,m-1}\notin I'$, then $-b+c\in I'$. Therefore $$b=c-(-b+c)\in I''+I'=I,$$
as required for the claim.

To conclude the proof of the statement, since $J:P$ is a prime ideal it is necessarily an associated prime of prime of $R/J$. As pointed out at the end of the proof of the previous item, $P\subset J:P$, hence $J:P$ is an embedded prime of $R/J$. Moreover, this also gives $P^2\subset J$, hence $J$ defines a double structure on the irreducible variety defined by $P$.

Let $\mathcal{Q}$ denotes the embedded component of $J$ with radical $J:P$ and let 
$\mathcal{Q}'$ denote the intersection of the remaining embedded components of $J$.
From $J=P\cap \mathcal{Q}\cap \mathcal{Q}'$ we get 
$$J:P=(\mathcal{Q}:P)\cap (\mathcal{Q}':P),$$
in particular, passing to radicals,  $J:P\subset \sqrt{\mathcal{Q}'}$.
This shows that $\mathcal{Q}$ is the unique embedded component of codimension $\leq 2(m-1)$, while the corresponding geometric component is supported on a linear subspace, as claimed in this item.

\medskip

{\rm (iv)} The ideal $P$ of submaximal minors specializes from the generic symmetric case, hence it is linearly presented as in that case.
In addition, the $k$-subalgebra generated by the submaximal minors has maximal dimension. Indeed, this is an immediate consequence of Theorem~\ref{cloning_sym} (iii) as the $k$-subalgebra generated by the partial derivatives of $f$ has dimension ${m+1\choose 2}-1$ (maximal possible)  and is a subalgebra thereof. 
(One could alternatively invoke \cite[Lemma 3.5]{LiZaSi}).

 Therefore, \cite[Theorem 3.2]{AHA} yields that the submaximal minors define a birational map onto the image.
In addition, since the minimal number of generators of $P$ is one plus the one of $J$, and the latter is generated by algebraically independent elements, it follows that the image of the birational map is a hypersurface. 

We set ourselves to describe the defining equation of this hypersurface.

Take the $m\times m$ generic symmetric matrix of the target coordinates  $(y_{i,j})_{1\leq i\leq j\leq m}$ and let $\mathbb{D}_{m-1,m-1}$ and $\mathbb{D}_{m,m}$ denote the respective cofactors of the $(m-1,m-1)$th and the $(m,m)$th entries of the target matrix.

Consider the matrix identity
		\begin{equation}\label{adj_eq}
		{\rm adj}({\rm adj}(\mathcal{SC}))=f^{m-2}\cdot\mathcal{SC},
		\end{equation}
where adj$(\_)$ denotes the transposed matrix of cofactors.

Looking at the right-hand side matrix of (\ref{adj_eq}) in more detail one sees that the entries in slots $(m-1,m-1)$ and $(m,m)$ are the same element, namely, $f^{m-2} x_{m-1,m-1}$.
Since the corresponding entries on the left-hand side matrix are
$\mathbb{D}_{m-1,m-1}(\underline{\Delta})$ and $\mathbb{D}_{m,m}(\underline{\Delta})$, respectively, one gets $(\mathbb{D}_{m,m}-\mathbb{D}_{m-1,m-1})(\underline{\Delta})=0$, as required.

But $\mathbb{D}_{m,m}-\mathbb{D}_{m-1,m-1}$ is clearly an irreducible polynomial in the target coordinate ring $k[y_{i,j}| 1\leq i\leq j\leq m]$. Therefore, we are through.
		
		\medskip
		
(v) As already observed in the proof of the previous item, $J$ is generated by the cofactors $\{\Delta_{l,h}|\,(l,h)\neq (m-1,m-1), (l,h)\neq (m,m)\}$
and the additional form $\Delta_{m,m}+\Delta_{m-1,m-1}$.

Now,  from (iv) the reduction number of a minimal reduction of $P$ is $m-2$. Thus, to conclude, it suffices to prove that $P^{m-1}\notin JP^{m-2}$.   
We will show that the element $\Delta_{m,m}^{m-1}$ of $P^{m-1}$ does not belong to  $JP^{m-2}$.
Supposing otherwise, we can write a polynomial relation of degree $m-1$ on the generators of $P$, namely, 
		
		\begin{equation} \label{reduction-relation1}
		\Delta_{m,m}^{m-1}=\sum_{1\leq l\leq h \leq m \atop (l,h)\neq (m-1,m-1),(m,m)} \Delta_{l,h}Q_{l,h}(\underline{\Delta})+(\Delta_{m-1,m-1}+\Delta_{m,m})Q(\underline{\Delta})
		\end{equation}
\noindent where $Q_{l,h}(\underline{\Delta})$ and $Q(\underline{\Delta})$ are homogeneous polynomial expressions of degree $m-2$ in the set $\underline{\Delta}=\{\Delta_{i,j}\,| \,1\leq i\leq j\leq m\}$ of the cofactors (generators of $P$).

		So the corresponding form of degree $m-1$ in $k[y_{i,j}| 1\leq i\leq j\leq m]$ is a scalar multiple of the polynomial $\boldsymbol H:=\mathbb{D}_{m,m}-\mathbb{D}_{m-1,m-1}$ obtained in the previous item.   We argue that this is impossible.

	Observe that the sum $$\sum_{(l,h)\neq (m-1,m-1)\atop (l,h)\neq (m,m)} \Delta_{l,h}Q_{l,h}(\underline{\Delta})$$  does not contain any nonzero terms of the form $\alpha\Delta_{m,m}^{m-1}$ or $\beta\Delta_{m-1,m-1}\Delta_{m,m}^{m-2}$. 	 Now, if these two terms appear in the second summand $(\Delta_{m-1,m-1}+\Delta_{m,m})Q(\underline{\Delta})$ they must have the same scalar coefficient, say, $c\in k$. 
		Bring the first of these to the left-hand side of (\ref{reduction-relation1}) to get a polynomial relation of $P$ having a term $(1-c)y_{m,m}^{m-1}$.  If $c\neq 1$, this is a contradiction  because any term of $\boldsymbol H$ has degree at most 1 in the  variable $y_{m,m}$.
			
		On the other hand, if $c=1$ then we still have a polynomial relation of $P$ having a term $y_{m-1,m-1}y_{m,m}^{m-2}$. Now, if $m>3$ this is again a contradiction due to the nature of $\boldsymbol H$ as the nonzero term of the latter has degree at most 1 in the variable $y_{m,m}$. 
		Finally, if $m=3$ a direct checking shows that the  monomial $y_{2,2}y_{3,3}$ can not to be the support of a nonzero term in $\boldsymbol H$. 
		This concludes the proof of the statement. 
		\qed
		
%		\begin{Remark}\rm
			%In the fully generic case \cite{LiZaSi} the primeness of the ideal of submaximal minors of the cloned matrix used a result of Eisenbud drawn upon the $2$-generic property.  Since the generic symmetric matrix is not $2$-generic, we were forced to devise an alternative to prove the difficult part of the assertion in (i) above.
	%	\end{Remark}

%\%\%\%\%\%\%\%\%\%\%\%\% %\%\%\%\%\%\%\%\%\%\%\%
\section{Sparsing with strategic zeros}\label{zeros}

In this part we fix integers $m,r$ with $1\leq r\leq m-2$ and consider the following degeneration of the $m\times m$ generic symmetric matrix:

{\small
	\begin{equation*}
	\left(
	\begin{array}{cccccccc}
	x_{1,1}&\ldots & x_{1,m-r}& x_{1,m-r+1} &  x_{1,m-r+2}&\ldots & x_{1,m-1} & x_{1,m}\\
	\vdots & \ldots & \vdots & \vdots & \vdots & \ldots & \vdots &\vdots \\
	x_{1,m-r}& \ldots & x_{m-r,m-r}& x_{m-r,m-r+1} & x_{m-r,m-r+2}&\ldots & x_{m-r,m-1} & x_{m-r,m}\\
	x_{1,m-r+1}&\ldots & x_{m-r,m-r+1}& x_{m-r+1,m-r+1} & x_{m-r+1,m-r+2} &\ldots & x_{m-r+1,m-1} & 0\\
	x_{1,m-r+2}&\ldots & x_{m-r,m-r+2}& x_{m-r+1,m-r+2} & x_{m-r+2,m-r+2} &\ldots & 0 & 0\\
	\vdots & \ldots & \vdots & \vdots & \vdots&\iddots &\vdots & \vdots\\
	x_{1,m-1} &\ldots & x_{m-r,m-1}&  x_{m-r+1,m-1}& 0&\ldots & 0 & 0\\[4pt]
	%\hline\\ [-7pt]
	x_{1,m} &\ldots & x_{m-r,m}& 0& 0&\ldots & 0 & 0\\ [3pt]
	\end{array}
	\right)
	\end{equation*}
}

Assuming that $m$ is fixed in the context, let us denote the above matrix by $\mathcal{S}(r)$.

\begin{Remark}\label{partials_are_cofactors}\rm
	By Proposition~\ref{GolMar}, any  cofactor $\Delta_{i,j}$ of $\mathcal{S}(r)$ such that $i+j\leq 2m-r$ is a scalar multiple of the partial derivative $\partial f/\partial x_{i,j}$ (the scalar is actually $\pm 1/2$).
\end{Remark}

Quite a bit of the arguments employed in the subsequent results work as well for the generic symmetric matrix (i.e., $r=0$). However, since this case is well-known in the literature, we will throughout assume that $r\geq 1$.

\begin{Lemma}\label{basic_count}
	Fix integers $m,r$ with $1\leq r\leq m-2$, and let $\mathcal S$ and $\mathcal{S}(r)$ denote, respectively, the generic $m\times m$ square matrix and its degeneration as above. Then the number of distinct entries in $\mathcal S$ set to zero to get $\mathcal{S}(r)$ is
	\begin{equation}\label{number_of_zeros}
	\mathfrak o (r):=
	\left\{\begin{array}{ccc} \frac{(r+1)^2}{4} & {\rm if} & r \;\text{\rm is odd}\\ [5pt]
	\frac{r(r+2)}{4} & {\rm if} &  r \;\text{\rm is even}
	\end{array}
	\right.
	\end{equation}
	In particular, the number of distinct variables appearing as entries in $\mathcal{S}(r)$ is ${{m+1}\choose 2}- \frac{(r+1)^2}{4}$ {\rm (}respectively, ${{m+1}\choose 2}- \frac{r(r+2)}{4}${\rm )} if $r$ is odd {\rm (}respectively, even{\rm )}.
\end{Lemma}
\demo
To count the distinct variables appearing as entries in the degeneration sector of $\mathcal S$ we may proceed, e.g., column-wise from right to left. 
If $r$ is even, we get the summation
$r+(r-2)+(r-4)+\cdots + 2$ as many such variables.
In other words,
$$\sum_{l=1}^{r/2}2l= 2 \sum_{l=1}^{r/2}l=2 {{r/2+1}\choose 2}=r(r+2)/4.$$

If $r$ is odd, get instead 
$$\sum_{l=1}^{(r+1)/2} (2l-1)=2 {{(r+1)/2+1}\choose 2}-(r+1)/2=(r+1)^2/4.\quad\quad\quad\quad $$%\square$$
\qed

\subsection{The structure of the polar image}

\begin{Proposition}\label{gradient-zeros-sym}
Let $R=k[\xx]$ denote the polynomial ring in the nonzero entries of $\mathcal{S}(r)$, with $1\leq r\leq m-2$ and set $f :=\det \mathcal{S}(r)$. 
One has:
\begin{enumerate}
\item[{\rm (a)}] $f$ is irreducible.
\item[{\rm (b)}] Let $J\subset R$ denote the gradient ideal of $f$. Then
 $$ {\rm cod}(J)= \left\{
	\begin{array}{ll}
	2 & \mbox{if $m-r=2$}\\
	3 & \mbox{otherwise.}
	\end{array}
	\right.
	$$ 
	In particular, $R/(f)$ is a normal domain if and only if $m-r\geq 3$.	
\end{enumerate}
\end{Proposition}
\demo (a) As a basic preliminary, we argue that the initial term of $f$ in the revlex monomial order is the product of the entries along the main anti-diagonal.
Indeed, this is true for the generic symmetric determinant and since the null entries of $\mathcal{S}(r)$ do not interfere in the anti-diagonal product, this remains true for $\mathcal{S}(r)$ as well -- and in fact, for any of its individual cofactor $\Delta_{i,j}$ such that $i+j\leq 2m-r$ (note that for the lex order this passage fails).

We will induct on $m$.
The initial step of the induction ($m=3$) will be subsumed in the general step -- alternatively, the determinant for this size is essentially the cubic polynomial with vanishing Hessian devised by Gordan and Noether as a counter-example to Hesse's famous mistaken assertion.

By the Laplace expansion along the first row, since $x_{1,1}$ only appears once and on the first row, one sees that $f=x_{1,1}f_1+g$, where $f_1$ is the determinant of the symmetric degeneration of order $(m-1)\times (m-1)$, for the same $r$, obtained by omitting the first row and the first column of the original $\mathcal{S}(r)$,
and both $f_1$ and $g$ belong to the proper $k$-subalgebra $R'\subset R$ omitting the variable $x_{1,1}$.

To show that $f$ is irreducible it suffices to prove that it is a primitive polynomial considered as a polynomial of degree $1$ in $R'[x_{1,1}]$.
In other words, it suffices to show that the ideal $(f_1,g)\subset R'$ has codimension $2$.

Now, if $r\leq m-3$ then $f_1$ is irreducible by the inductive hypothesis. On the other hand, since ${\rm in}(f)$ is a summand of $g$ as well, one must have ${\rm in}(g)={\rm in}(f)$. Since, by a similar token, ${\rm in}(f_1)$ is the product along the main anti-diagonal of the corresponding submatrix, it follows that $g$ cannot be a multiple of $f_1$.
This takes care of the case where $r\leq m-3$.
For the case where $r=m-2$, $f_1$ is itself a product along a main anti-diagonal, in which case it is clear that ${\rm in}(f)$ and ${\rm in}(f_1)$ have no common variables.

(b) Since $f$ is irreducible (radical would suffice), the gradient $J$ has codimension at least $2$.
On the other hand, $J$ is contained in the ideal of $(m-1)$-minors and the latter has codimension at most $3$, hence so does $J$. In the case $m-r=2$,  it is quite evident  that $J\subset \left(x_{1,m},x_{2,m} \right)$.
Therefore, one is done with this case.

%Consider the initial ideal ${\rm in}(J)$ of $J$ in the revlex monomial order.  We know that the initials terms of the submaximal minors of $\mathcal{S}$ are the products along the main anti-diagonal. This is not affected by the present degeneration.  By Proposition~\ref{GolMar}, the  cofactors $\Delta_{i,j}$ such that $i+j\leq 2m-r$ are scalar multiple of the partial derivatives. In particular $\Delta_{1,1}$ and $\Delta_{1,m}$ belongs to $J$. Thus, if $a_{i,j}$ denotes the $(i,j)$-entry  of $\mathcal{S}(r)$, then $\prod_{i+j=m+1}a_{i,j}$ and $\prod_{i+j=m+2}a_{i,j}$ belong to ${\rm in}(J)$. This ensures  that ${\rm in }(J)$  has codimension at least 2  and complete the case $m-r=2$.

Now, we assume  that $m-r\geq 3$. Let $a_{i,j}$ denote the $(i,j)$-entry  of $\mathcal{S}(r)$. We will show that $\mathfrak{p} \cdot \prod_{i+j=m}a_{i,j}$ belongs to ${\rm in}(J)$, for certain monomial $\mathfrak{p}$ involving in its support only variables $x_{i,j}$ such that $i+j=2m-r$.
Now, by Remark~\ref{partials_are_cofactors}, $\Delta_{1,1}$ and $\Delta_{1,m}$ belongs to $J$.
Therefore, by the argument at the beginning of the previous item, both  $\prod_{i+j=m+1}a_{i,j}$ and $\prod_{i+j=m+2}a_{i,j}$ belong to ${\rm in}(J)$ as well.
This then implies that the latter has codimension at least $3$.

 Claim: If $r$ is odd then $\mathfrak{p}=x_{\delta-1,\delta}^2\cdot x_{\delta -2,\delta+1}^2\cdots x_{m-r, m}^2\in J:I_{m-1}(\mathcal{S}(r))$, where $\delta =m-\frac{r-1}{2}$. 
 
 To see this, for each $k=0,\ldots,\frac{r-1}{2}$  consider the following  $(\delta +k)\times (\delta +k)$ submatrix of the cofactor matrix of $\mathcal{S}(r)$:
 
 $$\left( \begin{matrix}
 \Delta_{1,1} & \Delta_{1,2} & \ldots & \Delta_{1,\delta +k}\\
 \Delta_{2,1} & \Delta_{2,2} & \ldots &\Delta_{2,\delta +k} 	\\
 \vdots &  \vdots & \ldots & \vdots\\
 \Delta_{\delta +k,1} & \Delta_{\delta +k,2} &\ldots &\Delta_{\delta +k,\delta +k}
 \end{matrix}\right) $$
 
It will be shown that for any $k\in \{0,\ldots ,\frac{r-1}{2}\}$ one has $$x_{\delta-1,\delta}^2\cdot x_{\delta -2,\delta+1}^2\cdots x_{\delta -1-k,\delta+k}^2 \Delta_{i,j}\in J$$
 for all $1\leq i\leq j\leq \delta+k$  and then set $k= \frac{r-1}{2}$ to pull out the assertion in the claim.

 We induct on $k$. For the initial step, observe that there is one single cofactor not belonging to  $J$, namely, $\Delta_{\delta,\delta}$. 
 The cofactor formula yields the relation
 $$x_{1,\delta -1}\Delta_{\delta,1} +\cdots +x_{\delta-1,\delta-1}\Delta_{\delta,\delta-1} + x_{\delta-1,\delta}\Delta_{\delta,\delta}=0.$$
Since Remark~\ref{partials_are_cofactors} says that any cofactor $\Delta_{i,j}$ such that $i+j< 2\delta=2m-r+1$ belongs to $J$, then $x_{\delta-1,\delta}\Delta_{\delta,\delta} \in J$.
Therefore, certainly $x_{\delta-1,\delta}^2\Delta_{i,j} \in J$ for all $1\leq i\leq j\leq \delta$.
 
 For the inductive step, let $p_k=x_{\delta-1,\delta}^2\cdot x_{\delta -2,\delta+1}^2\cdots x_{\delta -1-k,\delta+k}^2$ and  suppose that  $p_k \cdot \Delta_{i,j}\in J$
for all $1\leq i\leq j\leq \delta+k$. This means that   $p_k \cdot \Delta_{i,j}\in J$
for all $1\leq i, j\leq \delta+k$ because $\Delta_{i,j}=\Delta_{j,i}$, by symmetry. Consider the following submatrix of the cofactor matrix of $\mathcal{S}(r)$: 
 $$\left( \begin{array}{cccc|c}
 \Delta_{1,1} & \Delta_{1,2} & \ldots & \Delta_{1,\delta +k} &  \Delta_{1,\delta +k+1}\\
 \Delta_{2,1} & \Delta_{2,2} & \ldots &\Delta_{2,\delta +k} &  \Delta_{2,\delta +k+1}	\\
 \vdots &  \vdots & \ldots & \vdots & \vdots\\
 \Delta_{\delta +k,1} & \Delta_{\delta +k,2} &\ldots &\Delta_{\delta +k,\delta +k}& \Delta_{\delta +k,\delta +k+1}\\ 
 \hline
 \Delta_{\delta +k+1,1} & \Delta_{\delta +k+1,2} &\ldots &\Delta_{\delta +k+1,\delta +k}& \Delta_{\delta +k+1,\delta +k+1} 
 \end{array}\right) $$
We will show that $p_{k+1}\Delta_{i,j}=p_k\cdot x_{\delta -k-2,\delta+(k+1)}^2 \Delta_{i,j}\in J$
for all $1\leq i\leq j\leq \delta+k+1$.  
\medskip

Since $p_k\cdot \Delta_{i,j}\in J$ for all $1\leq i\leq j\leq \delta+k$  and  $\Delta_{i,\delta +k+1} \in J$ for all $i\leq \delta-k-2$,  it remains to see that $$p_k\cdot x_{\delta-k-2,\delta+k+1}^2 \Delta_{i,\delta +k+1}\in J \; \mbox{for all} \; \delta-k-1\leq i\leq \delta+k+1.$$   

This time around, the cofactor formula yields the following relation:

\begin{equation}\label{onze} \sum_{j=1}^{\delta-k-2} x_{j,\delta -k-2}\Delta_{i,j} +\sum_{j=\delta-k-1}^{\delta+k+1}x_{\delta-k-2,j}\Delta_{i,j}=0,\; \mbox{for all}\; \delta-k-1\leq i\leq \delta+k+1.\end{equation}

For $i\neq \delta+k+1$ we multiply this equality by $p_k$ obtaining the following expression:

$$\sum_{j=1}^{\delta-k-2} p_k\cdot x_{j,\delta -k-2}\Delta_{i,j} +\sum_{j=\delta-k-1}^{\delta+k+1}p_k\cdot x_{\delta-k-2,j}\Delta_{i,j}=0,\; \mbox{for all}\; \delta-k-1\leq i\leq \delta+k.$$

Since, by hypothesis, $p_k\cdot \Delta_{i,j}\in J$ for all $1\leq i\leq j\leq  \delta+k$, this expressions give us $$p_k\cdot x_{\delta-k-2,\delta+k+1}\Delta_{i,\delta+k+1}\in J\; \mbox{for all} \;  \delta-k-1\leq i\leq \delta+k.$$

Using this, we conclude the statement for $i= \delta+k+1$. Indeed, we multiply the equality (\ref{onze}), with $i= \delta+k+1$,  by $p_k\cdot x_{\delta-k-2,\delta+k+1}$ obtaining that  
$$p_k\cdot x_{\delta-k-2,\delta+k+1}^2\Delta_{\delta+k+1,\delta+k+1}\in J.$$
\noindent This takes care of the claim.

In particular, it follows from the above claim that $\mathfrak{p} \cdot \Delta_{m,m}\in J.$ Since   
$${\rm in}(\mathfrak{p}\cdot \Delta_{m,m})=\mathfrak{p} \cdot \prod_{i+j=m}a_{i,j},$$ we are through.

\medskip

When $r$ is even  we consider $\mathfrak{p}=x_{\tilde{\delta}-1,\tilde{\delta}+1}^2\cdot x_{\tilde{\delta}-2,\tilde{\delta}+2}^2\cdots x_{m-r,m}^2$ instead, where  $\tilde{\delta}=m-\frac{r}{2}$. Proceeding  as above, we conclude that  $\mathfrak{p}\cdot \Delta_{m,m}\in J$ and therefore $J$ has codimension 3.
\qed

\bigskip

\begin{Remark}\rm
In contrast to the question in \cite[Conjecture 3.14]{LiZaSi} as to whether $R/J$ is Cohen--Macaulay  when $r=m-2$, here the picture is that $R/J$ has an embedded component of codimension $3$, thus taking us closer to the sub--Hankel environment (\cite{MAron}).
\end{Remark}

\begin{Theorem}\label{polar_structure_zeros}
	Keep the notation of {\rm Proposition~\ref{gradient-zeros-sym}}.
	The homogeneous coordinate ring  of the polar variety of $f$ in $\pp^{{{m+1}\choose 2}-\mathfrak o (r)-1}$ 
	is a symmetric ladder determinantal ring of dimension ${{m+1}\choose 2}-2\mathfrak o (r)$.  In particular, the analytic spread of $J$ is ${{m+1}\choose 2}-2\mathfrak o (r)$.
\end{Theorem}
\demo 
Let $\mathcal{L}=\mathcal{L}(m,r)$ denote  the set of boldface variables in the matrix $\widetilde{\mathcal{S}}(r)$ depicted below and let $I_{m-r}(\mathcal{L})$ stand for the ideal generated by the $(m-r)\times (m-r)$ minors of $\widetilde{\mathcal{S}}(r)$ involving only the variables in $\mathcal{L}$.  We  observe that $\widetilde{\mathcal{S}}(r)$ is of the form $\mathcal{S}(r)$ in the new variables $y_{i,j}$.

{\small
	\begin{equation*}%\label{generic-zeros}
%\widetilde{\mathcal{S}}(r)=
\left(
\begin{array}{cccccccc}
\mathbf{y_{1,1}}&\ldots & \mathbf{y_{1,m-r}}& \mathbf{y_{1,m-r+1}} &  \mathbf{y_{1,m-r+2}}&\ldots & \mathbf{y_{1,m-1}} & y_{1,m}\\
\vdots & \ldots & \vdots & \vdots & \vdots & \ldots & \vdots &\vdots \\
\mathbf{y_{1,m-r}}& \ldots & \mathbf{y_{m-r,m-r}}& \mathbf{y_{m-r,m-r+1}} &\mathbf{y_{m-r,m-r+2}}&\ldots & \mathbf{y_{m-r,m-1}} & y_{m-r,m}\\
\mathbf{y_{1,m-r+1}}&\ldots & \mathbf{y_{m-r,m-r+1}}& \mathbf{y_{m-r+1,m-r+1}} & \mathbf{y_{m-r+1,m-r+2}} &\ldots & y_{m-r+1,m-1} & 0 \\
\mathbf{y_{1,m-r+2}}&\ldots & \mathbf{y_{m-r,m-r+2}}& \mathbf{y_{m-r+1,m-r+2}} & y_{m-r+2,m-r+2} &\ldots &0  &0 \\
\vdots & \cdots & \vdots & \vdots &\vdots & \cdots &\vdots &\vdots \\
\mathbf{y_{1,m-1}} &\ldots & \mathbf{y_{m-r,m-1}}&  y_{m-r+1,m-1}&0 &\ldots &0  &0 \\[4pt]
%\hline\\ [-7pt]
y_{1,m} &\ldots & y_{m-r,m}& 0&0 &\ldots & 0 &0 \\ [3pt]
\end{array}
\right)
\end{equation*}}

%{\sc Lilly: I believe that if we take the polar variety as defined by a minimal set of generators of $J$ then we don't need the $1/2$'s in the matrix. I think that it suffices to eliminate the multiples of $2$ that appear as coefficients of the partial derivatives.}

\medskip

Since $\mathcal{L}$ can be extended to a $(m-1)\times (m-1)$ symmetric matrix of indeterminates, the ring $K[\mathcal{L}]/I_{m-r}(\mathcal{L}) $ is one of the so-called symmetric ladder determinantal  rings. 

%For further use, denote $\widehat{x_{i,j}}=x_{i,j}$ if $i\leq j$ and $\widehat{x_{i,j}}=x_{j,i}$ otherwise.
 
Let $\Delta_{j,i}$ stand for the (signed) cofactor of the $(i,j)$th entry of $\mathcal{S}(r)$ and let $f_{i,j}$ denote the $x_{i,j}$-derivative of $f$.
%  Observe that, by symmetry, $\Delta_{i,j}=\Delta_{j,i}$. Thus, by the  Proposition \ref{GolMar} we have  that the partial derivative $f_{i,i}$ of $f$ coincides with the (signed) cofactor de $x_{i,i}$ and  the $x_{i,j}$-derivative $f_{i,j}$ of $f$ coincides with the (signed) cofactor of  $x_{i,j}$ multiplied by 2, for all $i<j$.

\smallskip

{\sc Claim 1:} The homogeneous defining ideal of the image of the polar map of $f$ contains the ideal $I_{m-r}(\mathcal{L})$.

\smallskip

%Let $x_{i,j}$ denote a nonzero entry of $\mathcal{S}(r)$ 
Given  integers $1\leq i_1< i_2 <\ldots <i_{m-r}\leq m-1$, consider the following submatrix of the adjoint matrix of $\mathcal{S}(r)$:
{\large
	$$F= \left( \begin{matrix}
	\Delta_{1,i_1} &\widehat{\Delta}_{i_1,2}& \widehat{\Delta}_{i_1,3}&\cdots &\widehat{ \Delta}_{i_1,m-i_{m-r}+(m-r-1)}\\
	\Delta_{1,i_2} &\widehat{\Delta}_{i_2,2}& \widehat{\Delta}_{i_2,3}&\cdots & \widehat{\Delta}_{i_2,m-i_{m-r}+(m-r-1)}\\
	\vdots  &\vdots & \vdots &\cdots & \vdots \\
	\Delta_{1,i_{m-r}} &\widehat{\Delta}_{i_{m-r},2}& \widehat{\Delta}_{i_{m-r},3}&\cdots & \widehat{\Delta}_{i_{m-r},m-i_{m-r}+(m-r-1)}
	\end{matrix}\right).$$
}
Letting $$  C=\left(  \begin{matrix}
x_{1,i_{m-r}+1} & x_{1,i_{m-r}+2}& \cdots & x_{1,m-1} & x_{1,m}\\
\vdots & \vdots & \cdots & \vdots & \vdots \\
\widehat{x_{m-r,i_{m-r}+1}} & \widehat{x_{m-r,i_{m-r}+2}}& \cdots & x_{m-r,m-1} & x_{m-r,m}\\
\widehat{x_{m-r+1,i_{m-r}+1}} & \widehat{x_{m-r+1,i_{m-r}+2}}& \cdots & x_{m-r+1,m-1} &    0\\
\vdots & \vdots & \cdots &\vdots  &  \vdots\\
\widehat{x_{m-i_{m-r}+(m-r-2),i_{m-r}+1}} & \widehat{x_{m-i_{m-r}+(m-r-2),i_{m-r}+2}}& \cdots & 0 &    0\\
\widehat{x_{m-i_{m-r}+(m-r-1),i_{m-r}+1}} & 0 & \cdots & 0 &    0\\
\end{matrix}\right),$$
the cofactor identity ${\rm adj}( \mathcal{S}(r))\cdot \mathcal{S}(r) = \det(\mathcal{S}(r))\mathbb{I}_m$ yields the relation
$$F\cdot C=0.$$
Since the  columns of $C$ are linearly independent, it follows that the rank of $F$ is at most $m-i_{m-r}+(m-r-1) -(m-i_{m-r})=(m-r)-1$. In other words, the maximal minors of the matrix
$$Y=\left( \begin{matrix}
y_{1,i_1} &\widehat{y_{i_1,2}}& \widehat{y_{i_1,3}}&\cdots & \widehat{y_{i_1,m-i_{m-r}+(m-r-1)}}\\
y_{1,i_2} &\widehat{y_{i_2,2}}& \widehat{y_{i_2,3}}&\cdots & \widehat{y_{i_2,m-i_{m-r}+(m-r-1)}}\\
\vdots  &\vdots & \vdots &\cdots & \vdots \\
y_{1,i_{m-r}} &\widehat{y_{i_{m-r},2}}& \widehat{ y_{i_{m-r},3}}&\cdots & \widehat{y_{i_{m-r},m-i_{m-r}+(m-r-1)}}
\end{matrix}\right)$$
all vanish on the cofactors ${\bf\Delta}=\left\{\Delta_{i,j}\right\}$. Here,  $\widehat{y_{i,j}}=y_{i,j}$  if $i\leq j$ and    $\widehat{y_{i,j}}=y_{j,i}$ if $i\geq j$.

Since $f_{i,j}$ is a scalar multiple of $\Delta_{i,j}=\Delta_{j,i}$ (see  Proposition \ref{GolMar}), the maximal minors of $Y$  vanish on the partial derivatives of $f$, thus proving that the homogeneous defining ideal of the image of the polar map of $f$ contains the ideal $I_{m-r}(\mathcal{L})$.

{\sc Claim 2:} The codimension of  the ideal  $I_{m-r}(\mathcal{L}(m,r))$ is at least $\mathfrak o (r)$.

%For the proof of the claim, we can assume that $\widetilde{\mathcal{S}}(r)$ is of the form $\mathcal{S}(r)$ in the variables $y_{i,j}$.

%{\small
%	\begin{equation*}%\label{generic-zeros}
%	\left(
%	\begin{array}{cccccccc}
%	\textcolor{red}{y_{1,1}}&\ldots & \textcolor{red}{\frac{1}{2}y_{1,m-r}}& \textcolor{red}{\frac{1}{2}y_{1,m-r+1}} &  \textcolor{red}{\frac{1}{2}y_{1,m-r+2}}&\ldots & \textcolor{red}{\frac{1}{2}y_{1,m-1}} & \frac{1}{2}y_{1,m}\\
%	\vdots & \ldots & \vdots & \vdots & \vdots & \ldots & \vdots &\vdots \\
%	\textcolor{red}{y_{1,m-r}}& \ldots & \textcolor{red}{y_{m-r,m-r}}& \textcolor{red}{y_{m-r,m-r+1}} &\textcolor{red}{y_{m-r,m-r+2}}&\ldots & \textcolor{red}{y_{m-r,m-1}} & y_{m-r,m}\\
%	\textcolor{red}{y_{1,m-r+1}}&\ldots & \textcolor{red}{y_{m-r,m-r+1}}& \textcolor{red}{y_{m-r+1,m-r+1}} & \textcolor{red}{y_{m-r+1,m-r+2}} &\ldots & y_{m-r+1,m-1} & 0 \\
%	\textcolor{red}{y_{1,m-r+2}}&\ldots & \textcolor{red}{y_{m-r,m-r+2}}& \textcolor{red}{y_{m-r+1,m-r+2}} & y_{m-r+2,m-r+2} &\ldots &0  &0 \\
%	\vdots & \cdots & \vdots & \vdots &\vdots & \cdots &\vdots &\vdots \\
%	\textcolor{red}{y_{1,m-1}} &\ldots & \textcolor{red}{y_{m-r,m-1}}&  y_{m-r+1,m-1}&0 &\ldots &0  &0 \\[4pt]
%	%\hline\\ [-7pt]
%	y_{1,m} &\ldots & y_{m-r,m}& 0&0 &\ldots & 0 &0 \\ [3pt]
%	\end{array}
%	\right)
%	\end{equation*}} 

%Let us note that the codimension of this ladder ideal could be obtained by the general principle described in \cite{HeTr} (see also \cite{ladder2}), as done in the proof of Proposition~\ref{dim_dual}. However, in this structured situation we prefer to give an independent argument.

	We induct with the following inductive hypothesis: let $1\leq i\leq r-1;$ then for any $(m-i)\times (m-i)$ matrix of the form $\mathcal{S}(r)$, the ideal $I_{(m-i)-(r-i)}(\mathcal{L}(m-i,r-i))$ has codimension at least $\mathfrak o (r-i)$.  Note that $(m-i)-(r-i)=m-r,$ hence the size of the inner minor does not change in the inductive step. 
	
	We descend with regard to $i;$ thus, the induction step starts out at $i=r-1,$ hence  $m-i=m-r+1$ and we are in the situation of a $(m-r+1)\times (m-r+1)$ matrix of the form $\mathcal{S}(1)$.  Clearly, the ladder ideal $I_{(m-(r-1))-(r-(r-1))}{\mathcal{L}(m-r+1,r-(r-1))}$ is a principal ideal. Therefore, its codimension is $\mathfrak o (1)=1$ as desired. 
{\small
	\begin{equation*}%\label{generic-zeros}
	\left(
	\begin{array}{cccccccc}
	\mathbf{y_{1,1}}&\ldots & \mathbf{x_{1,m-r}}& \mathbf{y_{1,m-r+1}} &  \mathbf{y_{1,m-r+2}}&\ldots & \mathbf{y_{1,m-2}} & y_{1,m-1}\\
	\vdots & \ldots & \vdots & \vdots & \vdots & \ldots & \vdots &\vdots \\
	\mathbf{y_{1,m-r}}& \ldots & \mathbf{y_{m-r,m-r}}& \mathbf{y_{m-r,m-r+1}} &\mathbf{y_{m-r,m-r+2}}&\ldots & \mathbf{y_{m-r,m-2}} & y_{m-r,m-1}\\
	\mathbf{y_{1,m-r+1}}&\ldots & \mathbf{y_{m-r,m-r+1}}& \mathbf{y_{m-r+1,m-r+1}} & \mathbf{y_{m-r+1,m-r+2}} &\ldots & y_{m-r+1,m-2} & \\
	\mathbf{y_{1,m-r+2}}&\ldots & \mathbf{y_{m-r,m-r+2}}& \mathbf{y_{m-r+1,m-r+2}} & y_{m-r+2,m-r+2} & &  & \\
	\vdots & \ldots & \vdots & \vdots && & & \\
	\mathbf{y_{1,m-2}} &\ldots & \mathbf{y_{m-r,m-2}}&  y_{m-r+1,m-2}& & &  & \\[4pt]
	%\hline\\ [-7pt]
	y_{1,m-1} &\ldots & y_{m-r,m-1}& & & &  & \\ [3pt]
	\end{array}
	\right)
	\end{equation*}} 

To construct a suitable inductive predecessor, let $\widetilde{\mathcal{L}}$ denote the set of boldface variables above. Note that $\widetilde{\mathcal{L}}$ is of the form $\mathcal{L}(m-1,r-1)$ relative to a $(m-1)\times (m-1)$ matrix of the form $\mathcal{S}(r-1)$.  In particular $I_{m-r}(\widetilde{\mathcal{L}})$ is a Cohen--Macaulay prime ideal (see \cite{Conca} for primeness and Cohen--Macaulayness). By the inductive hypothesis, the codimension of $I_{m-r}(\widetilde{\mathcal{L}})$ is at least $\mathfrak{o}(r-1).$           %$(q_{r-1}+1)(q_{r-1}+e_{r-1}).$

Note that $\widetilde{\mathcal{L}}$ is a subset $\mathcal{L},$ hence there is a natural ring surjection:

$$S:=\frac{k[\mathcal{L}]}{I_{m-r}(\widetilde{\mathcal{L}})k[\mathcal{L}]}=\frac{k[\widetilde{\mathcal{L}}]}{I_{m-r}(\widetilde{\mathcal{L}})}[\mathcal{L}\setminus\widetilde{\mathcal{L}}]\surjects \frac{k[\mathcal{L}]}{I_{m-r}(\mathcal{L})}.$$

Since $\mathfrak{o}(r-1)+ \lceil\frac{r}{2}\rceil=\mathfrak{o}(r)$ it suffices to exhibit $\lceil\frac{r}{2}\rceil$ elements of $I_{m-r}(\mathcal{L})$ forming a regular sequence on the ring $S:=k[\mathcal{L}]/I_{m-r}(\widetilde{\mathcal{L}})k[\mathcal{L}].$

Consider the matrices 

\begin{equation}\label{Delta_i}
\left(\begin{array}{ccc|cccccccc}
x_{11}&\ldots&x_{1,m-r-1}&x_{1,m-i-1}\\
\vdots&\ddots&\vdots&\vdots\\
x_{1,m-r-1}&\ldots&x_{m-r-1,m-r-1}&x_{m-r-1,m-i-1}\\
\hline
x_{1,m-i}&\ldots&x_{m-r-1,m-i}&x_{m-r+i,m-i-1}
\end{array}\right)
\end{equation}
for $i=0,1,\ldots,\lceil\frac{r}{2}\rceil-1.$ Let $\Delta_i\in I_{m-r}(\mathcal{L})$ denote the determinant of the above matrix, for $i=0,1,\ldots,\lceil\frac{r}{2}\rceil-1.$

The claim is that $\boldsymbol\Delta=\{\Delta_0,\ldots,\Delta_{\lceil\frac{r}{2}\rceil-1}\}$ is a regular sequence on $S.$

Let $\delta$ denote the $(m-r-1)$-minor in the upper left corner of \eqref{Delta_i}. Clearly, $\delta$ is a regular element on $S$ as its defining is a prime ideal generated in degree $m-r.$ Therefore, it suffices to show that the localized sequence 

$$\boldsymbol\Delta_{\delta}=\{(\Delta_0)_{\delta},\ldots,(\Delta_{\lceil\frac{r}{2}\rceil-1})_{\delta}\}$$ is a regular sequence on $S_{\delta}.$ On the other hand, since $S$ is Cohen--Macaulay, it is suffices to show that $\dim S_{\delta}/\boldsymbol\Delta_{\delta}S_{\delta}=\dim S_{\delta}-\lceil\frac{r}{2}\rceil.$

Write $\XX=\{x_{m-r,m-1},x_{m-r+1,m-2},\ldots,x_{m-r+(\lceil\frac{r}{2}\rceil-1),m-(\lceil\frac{r}{2}\rceil-1)-1}\}.$ Note that, for every $i=0,\ldots,\lceil\frac{r}{2}\rceil-1,$ one has $(\Delta_i)_{\delta}=x_{m-r+i,m-i-1}+(1/\delta)\Gamma_i,$ with $x_{m-r+i,m-i-1}\in \XX$ and $\Gamma_i\in k[\mathcal{L}\setminus X].$ The association $x_{m-r+i,m-i-1}\mapsto -(1/\delta)\Gamma_i$ therefore defines a ring homomorphism  

$$k[\mathcal{L}]_{\delta}/(\boldsymbol\Delta_{\delta})=(k[\XX][\mathcal{L}\setminus\XX])_{\delta}/(\boldsymbol\Delta_{\delta})\simeq k[\mathcal{L}\setminus\XX]_{\delta}.$$

This entails a ring isomorphism 

$$\frac{S_{\delta}}{\boldsymbol\Delta_{\delta}S_{\delta}}\simeq\frac{k[\mathcal{L}\setminus X]_{\delta}}{(I_{m-r}(\widetilde{\mathcal{L}})k[\mathcal{L}\setminus \XX])_{\delta}}.$$

Thus, $\dim S_{\delta}/\boldsymbol\Delta_{\delta}S_{\delta}=\dim k[\mathcal{L}]_{\delta}-\lceil\frac{r}{2}\rceil-\codim I_{m-r}(\widetilde{\mathcal{L}})_{\delta}=\dim S_{\delta}-\lceil\frac{r}{2}\rceil.$

Therefore, $\codim(I_{m-r}(\mathcal{L}))$ is at least $\codim(I_{m-r}(\widetilde{\mathcal{L}}))+\lceil\frac{r}{2}\rceil=\mathfrak{o}(r).$

\bigskip

In order to show that $I_{m-r}(\mathcal{L}) $ is the homogeneous defining ideal  of the polar variety it suffices to show that the latter has codimension at most $\mathfrak{o}(r)$.
Since the dimension of the homogeneous coordinate ring of the polar variety coincides with the rank of the Hessian matrix of $f$, it now suffices to show
that the latter is at least $\dim R -\mathfrak{o}(r) = {{m+1}\choose {2}}- 2\mathfrak{o}(r)$.  

\medskip

Set $X:=\{ x_{i,j}\,|\, i+j=r+2,r+3,\ldots,2m-r \}$ and consider the set of partial derivatives  of $f$ with respect to the variables in $X$. Let $M$ denote the Jacobian matrix of these partial derivatives with respect to the variables in $X$. Observe that $M$ is a submatrix of size  $({{r+1}\choose {2}}-2\mathfrak{o})\times ({{r+1}\choose {2}}-2\mathfrak{o})$ of the Hessian matrix. We will show that $\det(M)\neq 0$.

Set  ${\bf v}:=\left\lbrace x_{i,j}\;|\; i+j=m+1\;\mbox{for}\; 1\leq i\leq j\leq m  \right\rbrace\subset X$, the set of variables along the main anti-diagonal of $\mathcal{S}(r)$. 

As already pointed out, the partial derivative of $f$ with respect to any $x_{i,i}$ coincides with the signed cofactor of $x_{i,i}$ and the partial derivative $f_{i,j}$ of $f$ with respect $x_{i,j}$, for $i<j$, is the (signed) cofactor of the variable in the $(i,j)$th entry multiplied by 2, i.e., $f_{i,j}=2\Delta_{i,j}$.  By expanding the cofactor of an entry in the set ${\bf v}$ one sees that there is a unique (nonzero) term whose support lives in ${\bf v}$ and the remaining terms have degree $\geq 2$ in the variables off ${\bf v}$. 
Similarly, the cofactor of a variable in ${X\setminus\bf v}$ has no term whose support lives in ${\bf v}$ and has exactly one (nonzero) term of degree 1 in the variables off $\vv$. In fact, if  $x_{i,j}\in X$ and $i+j\neq m+1$, one finds
\begin{eqnarray*}
	\Delta_{i,j} &=& \widehat{ x_{m+1-j,m+1-i}}\cdot \varpi\\
	& & \mbox{} + {\rm terms\; of\; degree\; at \;least\; 2\; off \;} \mathbf{v},
\end{eqnarray*}
where $\varpi$ is the product of the entries in the main anti-diagonal others than the entries in the slots $(i, m+1-i)$ and $(m+1-j,j)$. 

Consider the ring endomorphism $\varphi$ of $R$ that maps any variable in {\bf v}  to itself and  any variable off {\bf v} to zero.
By the preceding observation, applying $\varphi$ to any second partial derivative of $f$ involving only the variables of $X$ will return zero or a monomial supported on the variables in ${\bf v}$. 
Let $\widetilde{M}$ denote the matrix obtained by applying $\phi$ to the entries of  $M$.
Then any of its entries is either zero or a monomial supported on the variables in {\bf v}. 

We will  show that $\det(\widetilde{M}) $ is nonzero.
For this, consider the Jacobian matrix of the set of partial derivatives $\{f_v:v\in {\bf v}\}$ with respect to the variables in ${\bf v}$. Let $M_0$ denote the matrix obtained by applying $\phi$ to the entries of this Jacobian matrix by $\varphi$, viewed as a submatrix of $\widetilde{M}$. Up to permutation of rows and columns of $\widetilde{M}$, we may write

$$\widetilde{M}=
\left(
\begin{array}{cc}
M_0 & N_0 \\
N_1 & M_1
\end{array}
\right),
$$
for suitable $M_1$. Now, by the way in which the second partial derivatives of $f$ specialize via $\varphi$ as explained above, one must have $ N_0 = N_1 = 0$. Therefore, $\det(\widetilde{M}) = \det(M_0) \det(M_1)$, so it remains to prove the nonvanishing of these two subdeterminants. 
This is pretty much the same data as the end of the proof of Theorem~\ref{cloning_sym} (iii), so we conclude likewise.

This completes the proof of this item.
The supplementary assertion on the analytic spread of $J$ is clear since the dimension of the latter equals the dimension of the $k$-subalgebra generated by the partial derivatives.
\qed

\subsection{The structure of the submaximal minors}

%\textcolor{red}{\sc \% The idea is to prove two things: (1) the map defined by the cofactors is birational onto its image; (2) the image is a cone over the polar image. The other items of Theorem 3.7 of DG may be harder and sligtly different -- the preliminaries in DG are fairly involved. For example, Conj. 3.14 of DG is false here.}

In this subsection we consider aspects of the nature of the ideal generated by the submaximal minors of the degeneration $\mathcal{S}(r)$.
Some basic questions concerning its structure -- such as, e.g, its primeness as a function of $\mathfrak{o}$ -- are largely open.

\begin{Theorem}\label{submaximal_are_birational}
Let $I\subset R$ denote the ideal generated by the $(m-1)$-minors of $\mathcal{S}(r)$.
	
Then$:$
\begin{enumerate}
	 	\item[{\rm (b)}] $I$ is Cohen--Macaulay of codimension $3$.
		\item[{\rm (b)}] $I$ has maximal linear rank.
		\item[{\rm (c)}] The minors define a birational map $\pp^{{{m+1}\choose {2}}-\mathfrak{o}(r)-1}\dasharrow \pp^{{{m+1}\choose {2}}-1}$ onto a cone over the polar variety of $f$ with vertex cut by $\mathfrak{o}(r)$ coordinate  hyperplanes.
\end{enumerate}
\end{Theorem}

\demo (a) We claim that $I$ is a specialization from the generic symmetric case.
For this, one uses the same piece of proof as in \cite[Section 3.2]{LiZaSi}, using that the initial terms of the submaximal minors in the revlex order (respecting the rows) are the products of the variables along the corresponding anti-diagonals.
In addition, the proof that the distinct variables mapped to zero form a regular sequence on the ideal of submaximal minors of the generic symmetric matrix works just the same.

It follows immediately that $I$ is Cohen--Macaulay of codimension $3$.

(b) As in part (a), since $I$ is a specialization from the generic symmetric case it has the same free resolution data; in particular, it is even a linearly presented ideal.

An alternative argument for showing only maximality of the linear rank, without drawing on the Gr\"obner basis line of argument, goes as follows, pretty much along the same lines as the proof of Theorem~\ref{cloning_sym} (iv). 

First, get a hold of the notation $\widehat{x_{i,j}}$ from (\ref{defining_hat}). Let $\Delta_{i,j}$ stand  for the $(i,j)$-cofactor of $\mathcal{S}(r)$. The Cauchy cofactor formula  $$\mathcal{S}(r)\cdot {\rm adj}(\mathcal{S}(r))=\det(\mathcal{S}(r))\cdot \mathbb{I}_m$$
gives  the following linear relations involving the cofactors  of  $\mathcal{S}(r)$:

\begin{equation}
\label{eqsym1}
\left\{\begin{array}{ll}
\sum_{j=1}^m \widehat{x_{i,j}}\Delta_{j,1}=0& \kern-25pt\mbox{for $2\leq i \leq m-r$;}\\ [2pt]
\sum_{j=1}^{m-l} \widehat{x_{m-r+l,j}}\Delta_{j,1}=0 & \mbox{for  $1\leq l \leq r$;}
\end{array}
\right.
\end{equation}
\noindent with $m-1$ such relations;

\begin{equation}
\label{eqsym2}
\left\{\begin{array}{ll}
\sum_{j=1}^m \widehat{x_{i,j}}\Delta_{j,k}=0& \kern-25pt\mbox{for $2\leq k\leq m-r$ and  $k-1\leq i \leq m-r$ ($k\neq i$);}\\ [2pt]
\sum_{j=1}^{m-l} \widehat{x_{m-r+l,j}}\Delta_{j,k}=0& \mbox{for  $2\leq k\leq m-r$ and $1\leq l \leq r$;}
\end{array}
\right.
\end{equation}
\noindent with $(m-r-1)(m+r)/2$ such relations; 
\begin{equation}
\label{eqsym3}
\left\{\begin{array}{ll}
\sum_{j=1}^{m-k} \widehat{x_{m-r+k,j}}\Delta_{j,m-r+l}=0& \kern-70pt \mbox{for $1\leq l\leq r$ and  $0\leq k \leq r-l$ ($k\neq l$);}\\ [2pt]
\sum_{j=1}^{m-l} \widehat{x_{m-r+l,j}}\Delta_{j,m-r+l} - \sum_{j=1}^{m}\widehat{x_{1,j}}\Delta_{j,1}=0& \mbox{for  $1\leq l \leq \lfloor\frac{r}{2}\rfloor$;}
\end{array}
\right.
\end{equation}
\noindent with ${{r+1}\choose {2}}$ such relations.

Observing that $\Delta_{i,j}=\Delta_{j,i}$, we list the cofactors of the following form:
$$\Delta_{1,1}, \Delta_{2,1},\ldots, \Delta_{n,1}, \rightsquigarrow \Delta_{2,2}, \Delta_{3,2}, \ldots , \Delta_{n,2}\rightsquigarrow \ldots \rightsquigarrow \Delta_{m-1,m-1},\Delta_{m,m-1}\rightsquigarrow \Delta_{m,m}$$
\noindent With this ordering the above linear relations translate into linear syzygies of $I$ collected in the following block matrix:

{\scriptsize
	$$\mathcal{M}=\left(\begin{array}{ccccc|cccc}
	\varphi_1 &   \ldots &         &              &                 &     &        &  &                                 \\
	\mathbf{0}_{m-1}^{m-1}       &\varphi_2 &\ldots   &              &                 &     &     &   &                                     \\
	\vdots    &\vdots    & \ddots  & \vdots       &                 &     &        &    &                               \\
	\mathbf{0}^{m-1}_{r+2}        & \mathbf{0}_{r+2}^{m-1}        & \ldots  & \varphi_{m-r-1}&                 &  &  &        &                                  \\ 
	\mathbf{0}^{m-1}_{r+1}        & \mathbf{0}_{r+1}^{m-1}        & \ldots  & 0_{r+1}^{r+2} &\varphi_{m-r}&                 &  &  &                                        \\
	\hline\\ [-10pt]
	\mathbf{0}_{r}^{m-1} & \mathbf{0}_{r}^{m-1} & \ldots  &  \mathbf{0}_{r}^{r+2}& \mathbf{0}_{r}^{r+1}	     & \Phi_1   &             &                 &     \\
	\mathbf{0}_{r-1}^{m-1} & \mathbf{0}_{r-1}^{m-1}  &  \ldots & \mathbf{0}_{r-1}^{r+2}     &\mathbf{0}_{r-1}^{r+1}  &    \mathbf{0}_{r-1}^{r}&\Phi_2  &  &                      \\
	\vdots    & \vdots & \ldots  & \vdots       &  \vdots          &  \vdots        &\vdots &  \ddots\\
	\mathbf{0}_{1}^{m-1} & \mathbf{0}_{1}^{m-1}& \ldots  &  \mathbf{0}_{1}^{r+2}	    & \mathbf{0}_{1}^{r+1} &    \mathbf{0}_{1}^{r}       & \mathbf{0}_{1}^{r-1}   & \ldots &   \Phi_{r}
	\end{array}
	\right),$$}

where:

\begin{itemize}
	\item $\varphi_1$ is the matrix obtained  from $\mathcal{S}(r)$ by  omitting its first column;
	\item  For  $k=2,\ldots,m-r$, $\varphi_k$ is  the matrix of size $(m-k+1)\times (m-k+1)$ obtained  from $\mathcal{S}(r)$ by omitting the rows $1,\ldots, k-1$ and the  columns $1,\ldots,k-2,k$;
	\item   For $l=1,\ldots, r$,  $\Phi_l$ is the transpose of the  $(r-l+1))\times(r-l+1))$ minor obtained from the following submatrix  of $\mathcal{S}(r)$
	$$\Phi=\left(\begin{matrix}
	x_{m-r,m-r+1} & x_{m-r,m-r+2} & \ldots & x_{m-r,m}\\
	x_{m-r+1,m-r+1} & x_{m-r+1,m-r+2} & \ldots & 0\\
	\vdots & \vdots & \vdots & \vdots \\
	\widehat{x_{m-r+1,m-1}} & 0 & \ldots & 0 	
	\end{matrix}\right) $$
	by omitting the first $k$ rows and the last $k$  columns.		
	\item $\mathbf{0}_l^c$ denotes  an $l\times c $ block of zeros.
	
\end{itemize}

Counting through the sizes of the various blocks, one sees that this matrix is ${{m+1}\choose {2}}	\times ({{m+1}\choose {2}}-1)$. Omitting its first row obtains a square block-diagonal submatrix where each block has nonzero determinant. Thus, the linear rank of $P$ is maximal.

(c)
By part (a), $I$ has maximal linear rank.  On the other hand, its analytic spread  is maximal by \cite[Lemma 3.3]{LiZaSi}. Therefore, by  \cite[Theorem 3.2]{AHA} the submaximal minors define a birational map onto the image. It remains to argue that the image is a cone over the polar variety. 

To see this note that the homogeneous inclusion $T:=k[J_{m-1}]\subset T':=k[(I_{m-1})_{m-1}]$ of $k$-algebras which are domains, where  $I_{m-1}$ is minimally generated  by the  generators of $J$
and by $\mathfrak{o}(r)$ additional generators, say, $f_1,\ldots, f_{\mathfrak{o}(r)}$, that is, $T'=T[f_1,\ldots, f_{\mathfrak{o}(r)}]$. On the other hand, one has  $\dim T'={{m+1}\choose {2}}-\mathfrak{o}(r)$ and, by Theorem~\ref{polar_structure_zeros}, $\dim T= {{m+1}\choose {2}}-2\mathfrak{o}(r)$. Therefore,  ${\rm tr.deg}_{k(T)} k(T)(f_1,\ldots, f_{\mathfrak{o}(r)})=\dim T'-\dim T=\mathfrak{o}(r)$, where $k(T)$  denotes the field of fractions of $T$. This means that $f_1,\ldots, f_{\mathfrak{o}(r)}$ are algebraically independent over $k(T)$  and, a fortiori, over $T$. This shows that $T'$ is a polynomial ring over $T$ in $\mathfrak{o}(r)$ indeterminates. \qed
	
\section{The structure of the dual variety}\label{DUAL}

In this section we study the dual variety of the determinant of a degeneration of the $m\times m$ generic symmetric matrix.
It is well-known that in the generic symmetric case the dual variety is ideal theoretically defined by the $2$-minors, hence has dimension ${{m+1}\choose 2}-1-{m\choose 2}=m-1$.

Inspired by passages in \cite{Landsberg}, we ask whether the dual variety of the determinant of any ``reasonable'' degeneration of the $m\times m$ generic symmetric matrix has dimension $m-1$.

The difficulty is given a working interpretation of ``reasonable'' may become slightly clear in the subsequent development, in which we consider some special classes.

\subsection{Degenerations with dual of dimension $\geq m-1$}

Let $m\geq 3 $  be an integer. Throughout let $\mathcal{S}$ denote the $m\times m$ generic symmetric matrix as depicted in (\ref{symgeneric}).

\begin{Definition}\rm
A coordinate-like degeneration of $\mathcal{S}$ is $(m+2)$-{\em ladderlike} if the defining endomorphism fixes any variable $x_{i,j}$ such that $i+j\leq m+2$.
\end{Definition}

\begin{Theorem}\label{dual_ladderlike}
	Let $\mathcal{DS}$ denote any $(m+2)$-{ladderlike} degeneration of $\mathcal{S}$ and let $f=\det(\mathcal{DS})$. Then $\dim(V(f))^*\geq m-1.$
\end{Theorem}

\demo A formula due to Segre  (cf. \cite{Segre}, as revisited in \cite[Lemma 7.2.7]{Frusso}) tell us that 
\begin{equation}\label{Segre_formula}
\dim(V(f)^*)={\rm rank}\;  H(f)\; ({\rm mod}\; f)-2
\end{equation}

\noindent where $H(f)$ denotes  the Hessian matrix of $f$. Consequently, we have to prove that the latter has rank at least  $m+1$ modulo $f$.
To see this, let  $\Theta$ denote the Jacobian matrix of the following set of partial derivatives of $f$
$$\mathfrak{F}:=\left\lbrace \frac{\partial f}{\partial x_{2,m}},\frac{\partial f}{\partial x_{1,m}},\frac{\partial f}{\partial x_{1,m-1}},\ldots,\frac{\partial f}{\partial x_{1,2}},\frac{\partial f}{\partial x_{1,1}}\right\rbrace $$ 
with respect to the variables  $\xx=\{x_{1,1},x_{1,2},x_{1,3},\ldots,x_{1,m},x_{2,m}\}.$ Up to permutation of rows and columns, $\Theta$ is an $(m+1)\times (m+1)$  submatrix of $H(f).$ 
Our goal is to show that  $\det \Theta$ does not vanish modulo $f$.

As an easy case of Proposition~\ref{GolMar}, $\frac{\partial f}{\partial x_{1,1}}=\Delta_{1,1} $ while $\frac{\partial f}{\partial x_{i,j}}=2\Delta_{i,j} $ for every $x_{i,j}\in \xx\setminus\{x_{1,1}\}$. Then, its clear that the partial derivatives  $ \frac{\partial f}{\partial x_{1,m}},\frac{\partial f}{\partial x_{1,m-1}},\ldots,\frac{\partial f}{\partial x_{1,1}} $ do not depend on $x_{1,1}$. Consequently,  
\begin{equation}\label{state_of_derivatives}
\left\{\begin{array}{cc}
\frac{\partial^2 f}{\partial x_{1,1}\partial x_{1,j}}=\frac{\partial^2 f}{\partial x_{1,j}\partial x_{1,1}}=0, & \text{\rm for} \;j=1,\ldots, m\\ [5pt]
\frac{\partial^2 f}{\partial x_{1,j}\partial x_{k,l}}=\frac{\partial^2 f}{\partial x_{k,l}\partial x_{1,j}} & \text{\rm is independent of $x_{1,1},$ for $j=1,\ldots, m$ and $x_{k,l}\in\xx$} 
\end{array}
\right.
\end{equation}
 Therefore, $\Theta$ has the following shape:	
$$\Theta=\left(\begin{array}{c|ccc}
\frac{\partial^2 f}{\partial x_{1,1}\partial x_{2,m}}& L\\
\hline
{\bf 0}&\Omega
\end{array}\right)$$
\noindent	where ${\bf 0}$ is the $m\times 1$ zero matrix and  $\Omega$ is the $m\times m$ Jacobian matrix of the set $\mathfrak{F}\setminus\{ \frac{\partial f}{\partial x_{2,m}}\}$ with respect to the variables $\{x_{1,2},x_{1,3},\ldots,x_{1,m},x_{2,m}\}.$ 
Then 
$$\det\Theta=\frac{\partial^2 f}{\partial x_{1,1}\partial x_{2,m}}\cdot \det\Omega$$ 
is independent of $x_{11}$ by (\ref{state_of_derivatives}).  Since $f$ depends on $x_{1,1}$, it follows that $\det \Theta$ does not vanish modulo $f$ provided it does not vanish on the polynomial ring. 	

To see that $\det \Theta\neq 0$ it suffices to show that $\det \Omega\neq 0$ -- this is because $\Delta_{1,1}$ depends effectively on $x_{2,m}$ and hence the second derivative above is nonzero. 
For this consider the ring endomorphism of $R$ that maps any variable in the set $\vv= \left\lbrace x_{i,j}\; | i+j=m+2 \right\rbrace$  to itself and any variable off $\vv$ to zero. 
It suffices to prove that the specialized matrix $\Omega(\vv)$ has a non-vanishing determinant. 
But direct inspection shows that the latter is the diagonal matrix with entries the monomials  $2D/a_{m+2-j,j},$ 
with $2<j<m$ and $2D/x_{2,m}$, where $D=\prod_{i+j=m+2}a_{i,j}$, with $a_{i,j}$ denoting the $(i,j)$-entry of the matrix $\mathcal{DS}$.
\qed

\begin{Remark}\rm
	We don't know under what conditions the bound of Theorem~\ref{dual_ladderlike} remains valid for an $(m+1)$-ladder-like degeneration. The bound is still right in Example~\ref{bad_cloning}.
\end{Remark}

\subsection{The dual variety: $\mathcal{MD}$-cloning case}

We refer to the material and notation of Section~\ref{CLONING}.
Recall that $\mathcal{MD}$-cloning stands for cloning along the main diagonal of the $m\times m$ symmetric matrix $\mathcal{S}$, with $m\geq 3$.
 
\begin{Theorem}\label{dual_cloning}{\rm ($m\geq 3$)}
	Let $\mathcal{SC}$ denote an $\mathcal{MD}$-cloning of $\mathcal{S}$ and let $f=\det (\mathcal{SC})$. Then the dimension of the dual variety of $V(f)$ is $m-1$.
\end{Theorem}
\demo It is clear that $\mathcal{SC}$ is an $(m+2)$-ladderlike degeneration. Therefore, by Theorem~\ref{dual_ladderlike}, it suffices to show that  $\dim V(f)^*\leq  m-1$. 

Let $P\subset k[\yy]:=k[y_{i,j}\,|\, 1\leq i\leq j\leq m,\;(i,j)\neq (m,m)]$ denote  the homogeneous defining ideal of the dual variety $V(f)^*$ in its natural embedding, that is, one has an isomorphism of graded $k$-algebras induced by the assignment $y_{i,j}\mapsto \partial f/\partial x_{i,j}$:
\begin{equation*}%\label{algevra_of_dual}
k[\yy]/P\simeq T/(f)\cap T,
\end{equation*}
where $T=k[ \partial f/\partial x_{i,j}\,|\, 1\leq i\leq j\leq m,\;(i,j)\neq (m,m)]$.

Let us first effect a preliminary reduction, as follows: for any pair $(i,j)$ of indices, set
$$c_{i,j}=\left\{
\begin{array}{cc}
1 & \text{\rm if $i=j$}\\
2 & \text{\rm if $i\neq j$}
\end{array}
\right.
$$
Then Proposition~\ref{GolMar} gives $\partial f/\partial x_{ij}=c_{ij}\Delta_{i,j}$ for $(i,j)\neq (m-1,m-1)$, and $\partial f/\partial x_{m-1,m-1}= \Delta_{m-1,m-1}+\Delta_{m,m}$.

Now apply the harmless change of variables in the polynomial ring $k[\yy]$ that fixes $y_{m-1,m-1}$ and maps  $y_{ij}\mapsto c_{ij}^{-1}y_{ij}$, for $(i,j)\neq (m-1,m-1)$. Let $Q\subset k[\yy]$ denote the image of $P$ under this change.
Clearly, $Q$ still defines the dual variety of $V(f)$ as the image modulo $(f)$ of the linear system of the following cofactors: 
$$\{\Delta_{ij}\,|\,1\leq i\leq j\leq m, \,(i,j)\neq (m,m), (m-1,m-1)\}\cup \{\Delta_{m-1,m-1}+\Delta_{m,m}\}.$$

We claim that $Q$ contains the ideal  generated by the  $2\times 2$ minors of the following symmetric ladder matrix:

$$\mathcal{L}=\left(\begin{array}{cccccccc}
y_{1,1}&y_{1,2}&\ldots& y_{1,m-2}&y_{1,m}&y_{1,m-1}\\
y_{1,2}&y_{2,2}&\ldots& y_{2,m-2}&y_{2,m}&y_{2,m-1}\\
\vdots&\vdots&\ddots&\vdots&\vdots\\
y_{1,m-2}&y_{2,m-2}&\ldots& y_{m-2,m-2}&y_{m-2,m}&y_{m-2,m-1}\\
y_{1,m}&y_{2,m}&\ldots& y_{m-2,m}&y_{m-1,m}&\\
y_{1,m-1}&y_{2,m-1}&\ldots& y_{m-2,m-1}
\end{array}\right).$$

To see this, consider the following relation obtained  from the cofactor identity: 

\begin{equation}\label{id_cofatores}
{\rm adj}(\mathcal{SC})\cdot\mathcal{SC}\equiv 0 \,({\rm mod \,}f).
\end{equation}
For each  pair of indices $(i,j)$ such that  $1\leq i<j\leq m,$  let  $F_{ij}$ denote the  $2\times m$ submatrix of ${\rm adj}(\mathcal{SC})$ consisting of its  $i$-th and $j$-th rows. In addition, let $C$ stand for the $m\times (m-1)$ submatrix of $\mathcal{SC}$ consisting of its $m-1$ first columns . Then \eqref{id_cofatores} give us the relations

$$F_{ij}C\equiv 0\,(\mbox{mod\,}f),$$
for each $1\leq i\leq j\leq m.$  From this,  since  the rank of $C$ modulo $f$ is still $m-1,$ the rank of every $F_{i,j}$  is necessarily 1. This shows that every $2\times 2$ minor  of ${\rm adj}(\mathcal{SC})$ vanishes modulo $f.$ Therefore, each such minor involving only the cofactors in the set 
$$\boldsymbol\Delta:=\{\Delta_{ij}\,|1\leq i\leq j\leq m, \,(i,j)\neq (m,m), (m-1,m-1)\}$$ 
induces a $2\times 2$ minor of $\mathcal{L}$  vanishing on  $\boldsymbol\Delta$. On the other hand, by construction we obtain this way all the $2\times 2$ minors of $\mathcal{L}.$ This proves the claim.

As an ideal in the  polinomial ring  $A=k[I_1(\mathcal{L})_1],$ generated by the entries of $\mathcal{L},$ one has that $I_2(\mathcal{L})$ is a Cohen--Macaulay  prime  ideal  and its codimension is ${m\choose2}-2$ (\cite[Corollary 1.9, Theorem 1.13]{Conca}).  
 Clearly, its extension to the full polinomial ring $B=k[\yy]$ is still prime of codimension ${m\choose2}-2.$ Thus, $Q$ contains a prime ideal of  codimension ${m\choose2}-2$

Now consider the quadric  $ h:= y_{1,1}y_{m-1,m-1}-y_{1,m-1}^2-y_{1,m}^2.$ Since $I_2(\mathcal{L})$ is generated in degree $2$ and none of its generators involves $y_{m-1,m-1}$, then $h$ cannot be a $k$-linear combination of these generators. Thus, $h\notin I_2(\mathcal{L})$. But since $I_2(\mathcal{L})$ is a prime ideal it follows that $Q$ has codimension at least ${m\choose2}-1$. Therefore, $\dim V(f)^*= {m\choose2}-2-\codim Q \leq {m\choose2}-2-({m\choose2}-1)=m-1$, as was to be shown. \qed

\begin{Remark}\rm
The precise structure of the dual variety in the $\mathcal{MD}$-cloning case is not very clear. Its homogeneous defining ideal is not Cohen--Macaulay and contains many minimal generators which are quadric trinomials. Although there is computational evidence that this ideal is minimally generated in degree $2$, we know as yet no proof of this presumed fact.
\end{Remark}

We end this part by considering the question as to whether $f$ is a factor of its Hessian determinant $h(f)$ with multiplicity $\geq 1$. 
In distinction to this effective multiplicity, the {\em expected multiplicity} (according to Segre) is in this setup defined as ${{m+1}\choose 2}-2-\dim V(f)^*-1={{m}\choose 2}-2$, where $V(f)^*$ denotes the dual variety to the hypersurface $V(f)$ (see \cite{CRS}).

\begin{Proposition}\label{expected_multiplicity}{\rm ($m\geq 3$)}
	Let $\mathcal{SC}$ denote an $\mathcal{MD}$-cloning of $\mathcal{S}$ and let $f=\det (\mathcal{SC})$. Then $f$ is a factor of its Hessian determinant $h(f)$ with effective multiplicity equal to the expected multiplicity ${{m}\choose 2}-2$.
\end{Proposition}
\demo By the main result of Theorem~\ref{dual_cloning} and (\ref{Segre_formula}), one has 
$$\rk H(f) \pmod{f}=m-1+2=m+1.$$
But since $m\geq 3$ then $m+1<{{m+1}\choose 2}-1$, as one readily verifies, and hence, $h(f)\equiv 0 \pmod{f}$.
Therefore, $h(f)$ is a multiple of $f$ and is nonzero by Theorem~\ref{cloning_sym} (iii).
At the other end, the multiplicity of $f$ as a factor of $h(f)$ is at least the expected multiplicity ${{m}\choose 2}-2$ (cf. \cite[Section 2.1, p. 16]{CRS}).
Counting degrees we find that the degree of the residual factor $h(f)/f^{{{m}\choose 2}-2}$ is $2$. Since $deg(f)=m\geq 3$, we are through.
\qed

\subsection{The dual variety: the sparse case}

In this part we refer to the sparse-like degeneration $\mathcal{S}(r)$ of Section~\ref{zeros}.
In contrast to the structural content of Theorem~\ref{dual_cloning} here the dual variety of $\det \mathcal{S}(r)$ will actually be a ladder determinantal variety, not merely a subvariety of one such.

Besides, it will produce additional examples where the codimension of the dual variety of a determinantal hypersurface $f$ in its polar variety can be arbitrarily large when the Hessian determinant $h(f)$ vanishes.
The advantage of these new examples is that the ambient dimension can be smaller than in others previously known examples (such as in \cite{LiZaSi} and \cite{GRS}).

\begin{Theorem}\label{dual_zeros} Let $0\leq r\leq m-2$ and $f:=\det \mathcal{S}(r)$. Then:
	\begin{enumerate}	
		\item[{\rm (a)}] The dual variety $V(f)^*$ of $V(f)$ is a ladder determinantal variety of dimension $m-1$ defined by $2$-minors; in particular it is  arithmetically Cohen--Macaulay and its codimension  in the polar variety of $V(f)$ is ${m \choose 2}-\mathfrak{o}(r)$.
		\item[{\rm(b)}] $V(f)^*$  is arithmetically Gorenstein if and only if $r=m-2$.
	\end{enumerate}
\end{Theorem}
\demo (a) It is clear that $\mathcal{S}(r)$ is an $(m+2)$-ladderlike degeneration. Once again, by Theorem~\ref{dual_ladderlike}, it suffices to prove that  $\dim V(f)^*\leq  m-1$. 

Let $P\subset k[\yy]:=k[y_{i,j}\,|\, 2\leq i+j\leq 2m-r,\,i\leq j]$ denote the homogeneous defining ideal of the dual variety $V(f)^*$ in its natural embedding.
Proceeding as in the beginning of the proof of Theorem~\ref{dual_cloning}, we assume that $P$ has been harmlessly replaced by an ideal $Q\subset k[\yy]$ such that
\begin{equation*}%\label{algevra_of_dual}
k[\yy]/P\simeq {T}/(f)\cap {T},
\end{equation*}
where ${T}:=k[\Delta_{i,j}\,|\,2\leq i+j\leq 2m-r,\,i\leq j]$ is a $k$-subalgebra generated by cofactors.

Now, it suffices to prove that $Q$ has codimension at least ${m\choose 2} - \mathfrak{o}(r)$ because then
$\dim V(f)^*\leq {m+1\choose 2} - \mathfrak{o}(r)-1-({m\choose 2} - \mathfrak{o}(r))=m-1$, as required.

We claim that $Q$  contains the symmetric ladder determinantal ideal generated by the $2\times 2$ minors of the following matrix
{\small
	\begin{equation*}
	\left(
	\begin{array}{cccccccc}
	y_{1,1}&\ldots & y_{1,m-r}& y_{1,m-r+1} &  y_{1,m-r+2}&\ldots & y_{1,m-1} & y_{1,m}\\
	\vdots & \ldots & \vdots & \vdots & \vdots & \ldots & \vdots &\vdots \\
	y_{1,m-r}& \ldots & y_{m-r,m-r}& y_{m-r,m-r+1} & y_{m-r,m-r+2}&\ldots & y_{m-r,m-1} & y_{m-r,m}\\
	y_{1,m-r+1}&\ldots & y_{m-r,m-r+1}& y_{m-r+1,m-r+1} & y_{m-r+1,m-r+2} &\ldots & y_{m-r+1,m-1} & 0\\
	y_{1,m-r+2}&\ldots & y_{m-r,m-r+2}& y_{m-r+1,m-r+2} & y_{m-r+2,m-r+2} &\ldots & 0 & 0\\
	\vdots & \ldots & \vdots & \vdots & \vdots&\iddots &\vdots & \vdots\\
	y_{1,m-1} &\ldots & y_{m-r,m-1}&  y_{m-r+1,m-1}& 0&\ldots & 0 & 0\\[4pt]
	%\hline\\ [-7pt]
	y_{1,m} &\ldots & y_{m-r,m}& 0& 0&\ldots & 0 & 0\\ [3pt]
	\end{array}
	\right)
	\end{equation*}
}

This piece is the same as in the proof of Theorem~\ref{dual_cloning}.
We repeat the steps for the reader's convenience.

Consider the following relation afforded by the cofactor identity:
\begin{equation}\label{relation_mod_f}
{\rm adj}(\mathcal{S}(r))\cdot\mathcal{S}(r)\equiv 0\,({\rm mod\,}f).
\end{equation}
Further, for each pair of integers $i,j$ such that $1\leq i<j\leq m $ let $F_{ij}$ denote the $2\times m$ submatrix of ${\rm adj}(\mathcal{S}(r))$ consisting of the $i$th and $j$th rows.
In addition, let $C$ stand for the $m\times (m-1)$ submatrix of $\mathcal{S}(r)$  consisting of its $m-1$ leftmost columns.
Then (\ref{relation_mod_f}) gives the relations
$$F_{ij}C\equiv0\,({\rm mod\,}f),$$
for all $1\leq i<j\leq m $.
From this, since the rank of $C$ modulo $(f)$ is obviously still $m-1$, the one of every $F_{i,j}$ is necessarily $1$.
This shows that every $2\times 2$ minor of ${\rm adj}(\mathcal{S}(r))$ vanishes modulo $(f)$.
Therefore,  each such minor involving only cofactors in the set $$\boldsymbol\Delta:=\{\Delta_{i,j}\,|\,2\leq i+j\leq 2m-r,\,i\leq j\}$$
induces a $2\times 2$ minor of $\mathcal{L}$ vanishing  on $\boldsymbol\Delta$. On the other side, by construction this procedure yields all the $2\times 2$ minors of $\mathcal{L}$.
This proves the claim.

\smallskip

Now, since $I_2(\mathcal{L})$ is a symmetric ladder determinantal ideal on a suitable symmetric generic matrix it is a Cohen-Macaulay prime ideal (\cite[Theorem 1.13]{Conca}). 
Moreover, its codimension is ${m+1\choose 2} - \mathfrak{o}(r)-1-(m-1)={m\choose 2} - \mathfrak{o}(r)$ as follows from \cite[Corollary 1.9]{Conca}.

Thus, $Q$ has at least this codimension as well.
On the other hand, one also has
$${\rm cod} \,Q={\rm cod}\, P={m+1\choose 2} - \mathfrak{o}(r)-1-\dim (V^{\ast})\leq {m\choose 2} - \mathfrak{o}(r)$$
as follows from the above claim. 
Therefore, $I_2(\mathcal{L})\subset Q$ are prime ideals with the same codimension, and hence, $I_2(\mathcal{L})=Q.$ 
It is clear that the original homogeneous ideal $P$ is also generated by the $2$-minors of a suitable symmetric ladder.

Finally, the supplement about the codimension of $V^{\ast}$ in the polar variety follows immediately from Theorem~\ref{polar_structure_zeros}.

\medskip

(b) By the previous item, the homogeneous defining ideal of the dual variety is generated by the $2\times 2$ minors of a ladder matrix such as above, up to tagged nonzero coefficients. Observe that  the smallest square matrix containing all the entries of the latter is an $m\times m$ matrix of the form $\mathcal{S}(r)$. By \cite[Theorem (b), p. 120]{Conca2},  the ladder ideal is  Gorenstein if only if the inner corners of the ladder  have indices $(i,j)$ satisfying the equality $i+j=m+1$.
In the present case, the inner corners have indices satisfying the equation 
$$i+j=m-r+ (m-1)=m-r+1+(m-2)=\cdots = m-1 +(m-r)=2m-r-1.$$
Clearly this common value equals $m+1$ if and only if $r=m-2$.
\qed

\end{document}